\documentclass[10pt]{article}

\usepackage{amsmath}
\usepackage{amsthm}
\usepackage{amssymb}
\usepackage{amscd}
\usepackage{amsfonts}  %  ist neu
\usepackage{makeidx}
\usepackage{makeidx}     % allows index generation
\usepackage{graphicx}    % standard LaTeX graphics tool
\usepackage{multicol}    % use for the two-column index

\DeclareMathAlphabet\EuR{U}{eur}{m}{n}
\SetMathAlphabet\EuR{bold}{U}{eur}{b}{n}

\begin{document}

%%%%%%%%%%%%%%%%%%% Zaehler %%%%%%%%%%%%%%%%%%%%%%%%

%%\newtoks\theorembodyfont
%%\theorembodyfont{\slshape}
\newtheorem{theorem}{Theorem}[section]
\newtheorem{lemma}[theorem]{Lemma}
\newtheorem{proposition}[theorem]{Proposition}
\newtheorem{definition}[theorem]{Definition}
\newtheorem{example}[theorem]{Example}
\newtheorem{remark}[theorem]{Remark}
\newtheorem{corollary}[theorem]{Corollary}
\newtheorem{conjecture}[theorem]{Conjecture}
\newtheorem{problem}[theorem]{Problem}
{\catcode`@=11\global\let\c@equation=\c@theorem}
\renewcommand{\theequation}{\thetheorem}
% Hier werden Gleichungen und Theoreme zusammen gezaehlt.
%Soll ein anderer Zaehler statt theorem verwendet werden
%(entspr. dem \newtheorem-Befehl), muss 2-mal theorem
%durch diesen Zaehler ersetzt werden. (Die Zeilen entsprechen
% der Zaehlung von \newtheorem{equation}[theorem]).

\renewcommand{\theenumi}{\alph{enumi}}
\renewcommand{\labelenumi}{(\theenumi)}

\makeatletter
\renewcommand{\@seccntformat}[1]{\csname the#1\endcsname.\hspace{1em}}
\makeatother
\newcommand{\tit}[1]{\begin{bf} \begin{center} \begin{Large}
\section{#1}
\label{sec: #1}
\end{Large}\end{center}\end{bf}
\nopagebreak}

%%%%%%%%%%%%%%%%%%%%%%%%%%%%%% Diagrams %%%%%%%%%%%%%%%%%%%%%%%%%%%%%%

\newcommand{\squarematrix}[4]{\left( \begin{array}{cc} #1 & #2 \\ #3 &
#4
\end{array} \right)}
\newcommand{\smallmx}[4]{\mbox{\begin{scriptsize}$\squarematrix{#1}{#2}
        {#3}{#4}$\end{scriptsize}}}

\newcommand{\comsquare}[8]{
\begin{center}
$\begin{CD}
#1 @>#2>> #3\\
@V{#4}VV @VV{#5}V\\
#6 @>>#7> #8
\end{CD}$
\end{center}}

%%%%%%%%%%%%%%%%%% Categories %%%%%%%%%%%%%%%%%%%%%

\let\sect=\S
\newcommand{\curs}{\EuR}
\newcommand{\CHAINCOMPLEXES}{\curs{CHCOM}}
\newcommand{\GROUPOIDS}{\curs{GROUPOIDS}}
\newcommand{\PAIRS}{\curs{PAIRS}}
\newcommand{\FGINJ}{\curs{FGINJ}}
\newcommand{\Or}{\curs{Or}}
\newcommand{\SPACES}{\curs{SPACES}}
\newcommand{\SPECTRA}{\curs{SPECTRA}}
\newcommand{\Sub}{\curs{Sub}}

%%%%%%%%%%%%%%% mathbb %%%%%%%%%%%%%%%%%%

\newcommand{\bbC}{{\mathbb C}}
\newcommand{\bbF}{{\mathbb F}}
\newcommand{\bbH}{{\mathbb H}}
\newcommand{\bbI}{{\mathbb I}}
\newcommand{\bbK}{{\mathbb K}}
\newcommand{\bbKO}{\mathbb{KO}}
\newcommand{\bbN}{{\mathbb N}}
\newcommand{\bbP}{{\mathbb P}}
\newcommand{\bbQ}{{\mathbb Q}}
\newcommand{\bbR}{{\mathbb R}}
\newcommand{\bbZ}{{\mathbb Z}}

\newcommand{\cala}{{\mathcal A}}
\newcommand{\calbh}{{\mathcal B}{\mathcal H}}
\newcommand{\calc}{{\mathcal C}}
\newcommand{\cald}{{\mathcal D}}
\newcommand{\cale}{{\mathcal E}}
\newcommand{\calf}{{\mathcal F}}
\newcommand{\calg}{{\mathcal G}}
\newcommand{\calh}{{\mathcal H}}
\newcommand{\calk}{{\mathcal K}}
\newcommand{\calp}{{\mathcal P}}
%%%%%%%%%%%%%%% mathbf %%%%%%%%%%%%%%%%%%

\newcommand{\bfA}{\ensuremath{\mathbf{A}}}
\newcommand{\bfa}{\ensuremath{\mathbf{a}}}
\newcommand{\bfb}{\ensuremath{\mathbf{b}}}
\newcommand{\bfCTR}{\ensuremath{\mathbf{CTR}}}
\newcommand{\bfDtr}{\ensuremath{\mathbf{Dtr}}}
\newcommand{\bfE}{\ensuremath{\mathbf{E}}}
\newcommand{\bff}{\ensuremath{\mathbf{f}}}
\newcommand{\bfF}{\ensuremath{\mathbf{F}}}
\newcommand{\bfg}{\ensuremath{\mathbf{g}}}
\newcommand{\bfHH}{\ensuremath{\mathbf{HH}}}
\newcommand{\bfHPB}{\ensuremath{\mathbf{HPB}}}
\newcommand{\bfI}{\ensuremath{\mathbf{I}}}
\newcommand{\bfi}{\ensuremath{\mathbf{i}}}
\newcommand{\bfK}{\ensuremath{\mathbf{K}}}
\newcommand{\bfL}{\ensuremath{\mathbf{L}}}
\newcommand{\bfr}{\ensuremath{\mathbf{r}}}
\newcommand{\bfs}{\ensuremath{\mathbf{s}}}
\newcommand{\bfS}{\ensuremath{\mathbf{S}}}
\newcommand{\bft}{\ensuremath{\mathbf{t}}}
\newcommand{\bfT}{\ensuremath{\mathbf{T}}}
\newcommand{\bfTC}{\ensuremath{\mathbf{TC}}}
\newcommand{\bfU}{\ensuremath{\mathbf{U}}}
\newcommand{\bfu}{\ensuremath{\mathbf{u}}}
\newcommand{\bfv}{\ensuremath{\mathbf{v}}}
\newcommand{\bfw}{\ensuremath{\mathbf{w}}}

\newcommand{\Rat}{\ensuremath{\mathbf{R}}\ensuremath{\mathbf{a}}\ensuremath{\mathbf{t}}}
%%%%%%%%%%%%%%%%%%%%% Operatornames %%%%%%%%%%%%%%%%%%%%%%%%%

\newcommand{\aut}{\operatorname{aut}}
\newcommand{\Bor}{\operatorname{Bor}}
\newcommand{\ch}{\operatorname{ch}}
\newcommand{\class}{\operatorname{class}}
\newcommand{\cok}{\operatorname{coker}}
\newcommand{\cone}{\operatorname{cone}}
\newcommand{\con}{\operatorname{con}}
\newcommand{\conhom}{\operatorname{conhom}}
\newcommand{\cyclic}{\operatorname{cyclic}}
\newcommand{\dirlim}{\operatorname{colim}}
\newcommand{\ev}{\operatorname{ev}}
\newcommand{\ext}{\operatorname{ext}}
\newcommand{\Gen}{\operatorname{Gen}}
\newcommand{\hur}{\operatorname{hur}}
\newcommand{\im}{\operatorname{im}}
\newcommand{\inj}{\operatorname{inj}}
\newcommand{\id}{\operatorname{id}}
\newcommand{\infl}{\operatorname{Infl}}
\newcommand{\ind}{\operatorname{ind}}
\newcommand{\Inn}{\operatorname{Inn}}
\newcommand{\Irr}{\operatorname{Irr}}
\newcommand{\Is}{\operatorname{Is}}
\newcommand{\map}{\operatorname{map}}
\newcommand{\MOD}{\operatorname{MOD}}
\newcommand{\mor}{\operatorname{mor}}
\newcommand{\Ob}{\operatorname{Ob}}
\newcommand{\op}{\operatorname{op}}
\newcommand{\pr}{\operatorname{pr}}
\newcommand{\point}{\operatorname{pt.}}
\newcommand{\res}{\operatorname{res}}
\newcommand{\topo}{\operatorname{top}}
\newcommand{\tors}{\operatorname{tors}}

%%%%%%%%%%%%%%%%%%% Limits %%%%%%%%%%%%%%%%%%%%%

\newcommand{\higherlim}[3]{\underleftarrow{\lim}_{#1}^{#2}#3}
\newcommand{\highercolim}[3]{\underrightarrow{\dirlim}_{#1}^{#2}#3}
\newcommand{\invlim}[2]{\higherlim{#1}{}{#2}}
\newcommand{\colim}[2]{\highercolim{#1}{}{#2}}

%%%%%%%%%%%%%%%%%% Families of subgroups %%%%%%%%%%%%%

\newcommand{\calfin}{{\mathcal F}\!{\mathcal I}\!{\mathcal N}}
\newcommand{\calmfin}{{\mathcal M}\!{\mathcal F}\!{\mathcal I}\!{\mathcal N}}

%%%%%%%%%%%%%%%%%%%%%%%%%% others %%%%%%%%%%%%%%%%%%%%%%%%%%%%%%%%%%%%%%

\newcommand{\pt}{\{\bullet\}}

\newcounter{commentcounter}
\newcommand{\comment}[1]                      %comment of the author
{\stepcounter{commentcounter}
{\bf Comment \arabic{commentcounter}}: {\ttfamily #1} }

\newcommand{\tabtit}[1]
{
\ref{sec: #1}. & #1
}

%%%%%%%%%%%%%%%%%%%%%%%%%  Begin of the text  %%%%%%%%%%%%%%%

\title{Rational Computations of the Topological {$K$}-Theory of
  Classifying Spaces of Discrete Groups}
\author{
Wolfgang L\"uck\thanks{\noindent email:
lueck@math.uni-muenster.de\protect\\
www: ~http://www.math.uni-muenster.de/u/lueck/\protect\\
FAX: 49 251 8338370\protect}\\
Fachbereich Mathematik\\ Universit\"at M\"unster\\
Einsteinstr.~62\\ 48149 M\"unster\\Germany}
\maketitle

%%%%%%%%%%%%%%%%%%%%%%%%%%%% Abstract  %%%%%%%%%%%%%%%%%%%%%%%%%%%%%%%%%%%
\typeout{-----------------------  Abstract  ------------------------}
\begin{abstract}
We compute rationally the topological (complex) K-theory of the classifying space $BG$
of a discrete group provided that $G$ has a cocompact $G$-$CW$-model for
its classifying space for proper $G$-actions.
For instance word-hyperbolic groups
and cocompact discrete subgroups of connected Lie groups satisfy this assumption.
The answer is given in terms of the group cohomology of $G$
and of the centralizers of finite cyclic subgroups of prime power order.
We also analyze the multiplicative structure.

\smallskip

\noindent
Key words: topological $K$-theory, classifying spaces of groups.\\
Mathematics Subject Classification 2000: 55N15.
\end{abstract}

%%%%%%%%%%%%%%%%%%%%%%%%%%%%%Introduction %%%%%%%%%%%%%%%%%%%%%%%%%%%%%%%%
\typeout{-----------------------  Introduction ------------------------}

\setcounter{section}{-1}
\tit{Introduction and Statements of Results}

The main result of this paper is:

\begin{theorem}[Main result] \label{the: main theorem}
Let $G$ be a discrete group. Denote by $K^*(BG)$
the topological (complex) K-theory of its classifying space $BG$.
Suppose that there is a
cocompact  $G$-$CW$-model for the classifying space $\underline{E}G$ for
proper $G$-actions.

Then there is a $\bbQ$-isomorphism
\begin{multline*}
\overline{\ch}^n_G \colon K^n(BG) \otimes_{\bbZ} \bbQ ~  \xrightarrow{\cong}
\\
\left(\prod_{i \in \bbZ} H^{2i+n}(BG;\bbQ)\right) \times
\prod_{p \text{ prime}} ~ \prod_{(g) \in \con_p(G)}
\left(\prod_{i \in \bbZ} H^{2i+n}(BC_G\langle g \rangle;\bbQ\widehat{_p})\right),
\end{multline*}
where $\con_p(G)$ is the set of conjugacy classes (g) of elements $g
\in G$ of order $p^d$ for some integer $d \ge 1$ and $C_G\langle g
\rangle$ is the centralizer of
the cyclic subgroup $\langle g \rangle$ generated by $g$.
\end{theorem}

The \emph{classifying space $\underline{E}G$ for proper $G$-actions}
is a proper $G$-$CW$-complex such that the $H$-fixed point set is
contractible for every finite subgroup $H \subseteq G$. It has the
universal property that for every proper $G$-$CW$-complex $X$ there is
up to $G$-homotopy precisely one $G$-map $f \colon X \to
\underline{E}G$. Recall that a $G$-$CW$-complex is proper if and only if
all its isotropy groups are finite, and is finite if and only if it is cocompact.
The assumption in Theorem~\ref{the: main theorem}
that there is a cocompact $G$-$CW$-model for the classifying space $\underline{E}G$ for
proper $G$-actions is satisfied
for instance if $G$ is word-hyperbolic in the sense of Gromov,
if $G$ is a cocompact subgroup of a Lie group with finitely many path
components, if $G$ is a finitely generated one-relator group,
if $G$ is an arithmetic group, a mapping class group of a compact surface or the group
of outer automorphisms of a finitely generated free group. For more information about
$\underline{E}G$ we refer for instance to~\cite{Baum-Connes-Higson(1994)} and~\cite{Lueck(2004a)}.
We will prove Theorem~\ref{the: main theorem} in Section~\ref{sec: Proof of the Main Result}.

We will also investigate the multiplicative structure on
$K^n(BG) \otimes_{\bbZ} \bbQ$ in Section~\ref{sec: Multiplicative Structures}.
If one is willing to complexify, one can show:

\begin{theorem}[Multiplicative structure] \label{the: Multiplicative structure}
Let $G$ be a discrete group. Suppose that there is a
cocompact  $G$-$CW$-model for the classifying space $\underline{E}G$ for
proper $G$-actions.

Then there is a $\bbC$-isomorphism
\begin{multline*}
\overline{\ch}^n_{G,\bbC} \colon K^n(BG) \otimes_{\bbZ} \bbC ~ \xrightarrow{\cong}
\\
\left(\prod_{i \in \bbZ} H^{2i+n}(BG;\bbC)\right) \times
\prod_{p \text{ prime}} ~ \prod_{(g) \in \con_p(G)}
\left(\prod_{i \in \bbZ} H^{2i+n}(BC_G\langle g \rangle;\bbQ\widehat{_p} \otimes_{\bbQ} \bbC)\right),
\end{multline*}
which is compatible with the standard multiplicative structure on
$K^*(BG)$ and the one on the target given by
$$
\left(a, u_{p,(g)}\right) \cdot \left(b, v_{p,(g)}\right)
~ = ~
\left(a\cdot b, (a \cdot v_{p,(g)} + b \cdot u_{p,(g)} + u_{p,(g)} \cdot v_{p,(g)})\right)
$$
for
\begin{eqnarray*}
(g) & \in & \con_p(G);
\\
a,b  & \in & \prod_{i \in \bbZ} H^{2i+*}(BG;\bbC);
\\
u_{p,(g)}, v_{p,(g)} & \in & \prod_{i \in \bbZ} H^{2i+*}(BC_G\langle g \rangle;\bbQ\widehat{_p} \otimes_{\bbQ} \bbC),
\end{eqnarray*}
and the structures of a graded commutative ring on  $\prod_{i \in \bbZ} H^{2i+*}(BG;\bbC)$ and
$\prod_{i \in \bbZ} H^{2i+*}(BC_G\langle g \rangle;\bbQ\widehat{_p} \otimes_{\bbQ} \bbC)$ coming from the
cup-product and the obvious $\prod_{i \in \bbZ} H^{2i+*}(BG;\bbC)$-module structure
on $\prod_{i \in \bbZ} H^{2i+*}(BC_G\langle g \rangle;\bbQ\widehat{_p} \otimes_{\bbQ} \bbC)$ coming from the
canonical  maps $BC_G\langle g \rangle \to BG$ and $\bbC \to \bbQ\widehat{_p} \otimes_{\bbQ} \bbC$.
\end{theorem}

In Section~\ref{sec: Weakening the Finiteness Conditions} we will
prove Theorem~\ref{the: main theorem} and Theorem~\ref{the: Multiplicative structure}
under weaker finiteness assumptions than stated above.

If $G$ is finite, we get the following integral computation of $K^*(BG)$.
Throughout the paper $R(G)$ will be the complex representation ring and
$\bbI_G$ be its augmentation ideal, i.e. the kernel of the ring homomorphism
$R(G) \to \bbZ$ sending $[V]$ to $\dim_{\bbC}(V)$.  If $G_p \subseteq G$ is a $p$-Sylow
subgroup, restriction defines a map $\bbI(G) \to \bbI(G_p)$. Let $\bbI_p(G)$ be the quotient
of $\bbI(G)$ by the kernel of this map. This is independent of the choice of
the $p$-Sylow subgroup since two $p$-Sylow subgroups of $G$ are conjugate.
There is an obvious isomorphism from $\bbI_p(G) \xrightarrow{\cong} \im(\bbI(G) \to
\bbI(G_p))$.
We will prove in Section~\ref{sec: The K-Theory of the Classifying Space of a Finite Group}

\begin{theorem}{\bf ($K$-theory of $BG$ for finite groups $G$).}
 \label{the: computation of K_*(BG) and K^*(BG) for finite G}
Let $G$ be a finite group.  For a prime $p$
denote by $r(p) = |\con_p(G)|$ the number of conjugacy classes $(g)$
of elements $g \in G$ whose order $|g|$ is $p^d$ for some integer $d \ge 1$.
Then there are isomorphisms of abelian groups
\begin{eqnarray*}
K^0(BG) & \cong & \bbZ \times \prod_{p\text{ prime}}  \bbI_p(G)
\otimes_{\bbZ} \bbZ\widehat{_p}
~ \cong ~
\bbZ \times \prod_{p\text{ prime}} (\bbZ\widehat{_p})^{r(p)};
\\
K^1(BG) & \cong & 0.
\end{eqnarray*}
The isomorphism $K^0(BG) \xrightarrow{\cong} \bbZ \times \prod_{p\text{ prime}}  \bbI_p(G)\otimes_{\bbZ} \bbZ\widehat{_p}$
is compatible with the standard ring structure on the source and the ring structure on the target given by
$$\left(m,u_p \otimes a_p\right) \cdot \left(n,v_p \otimes b_p\right) ~ = ~
\left(mn,(mv_p \otimes b_p + nu_p \otimes a_p + (u_pv_p) \otimes (a_pb_p)\right)
$$
for $m,n \in \bbZ$, $u_p,v_p \in \bbI_p(G)$ and $a_p, b_p \in \bbZ\widehat{_p}$
and the obvious multiplication in $\bbZ$, $\bbI_p(G)$ and $\bbZ\widehat{_p}$.
\end{theorem}

The additive version of Theorem~\ref{the: computation of K_*(BG) and K^*(BG) for finite G}
has already been explained in \cite[page~125]{Jackowski-Oliver(1996)}.
Inspecting \cite[Theorem~2.2]{Jackowski(1978)} one can also derive the ring structure. In
\cite{Kuhn(1987)} the $K$-theory of $BG$ with coefficients in the field $\bfF_p$ of
$p$ elements has been determined including the multiplicative structure.
The proof of Theorem~\ref{the: computation of K_*(BG) and K^*(BG) for finite G}
we will present here is based on the ideas of this paper.
We will and need to show a stronger statement about the pro-group
$\{\bbI_G/(\bbI_G)^{n+1}\}$ in Theorem~\ref{the: computation of AI_G/(AI_G)^{n+1}}
\eqref{the: computation of AI_G/(AI_G)^{n+1}: computing I/I^n}.

A version of Theorem~\ref{the: main theorem} for topological $K$-theory with coefficients
in the $p$-adic integers has been proved
by Adem~\cite{Adem(1992)}, \cite{Adem(1993b)} using the Atiyah-Segal completion theorem
for the finite group $G/G'$
provided that $G$ contains a torsionfree subgroup $G'$ of finite index.
Our methods allow to drop this condition, to deal with $K^*(BG) \otimes_{\bbZ} \bbQ$ directly and study
systematically the multiplicative structure for $K^*(BG) \otimes_{\bbZ} \bbC$. They are
based on the equivariant cohomological Chern character of~\cite{Lueck(2004i)}.

For integral computations of the $K$-theory and $K$-homology of
classifying spaces of groups we refer to~\cite{Joachim-Lueck(2005)}.

The paper is organized as follows:

\begin{tabular}{ll}
\tabtit{Borel Cohomology and Rationalization}
\\
\tabtit{Some Preliminaries about Pro-Modules}
\\
\tabtit{The K-Theory of the Classifying Space of a Finite Group}
\\
\tabtit{Proof of the Main Result}
\\
\tabtit{Multiplicative Structures}
\\
\tabtit{Weakening the Finiteness Conditions}
\\
\tabtit{Examples and Further Remarks}
\\
 & References
\end{tabular}
\bigskip

The author wants to the thank the Max Planck Institute for Mathematics in Bonn for its hospitality
during his stay from April 2005 until July 2005 when this paper was written.

%%%%%%%%%%%%%%%%%%%%%%%%%%%%%%%%%%%%%%%%%%%%%%%%%%%%%%%%%%%%%%%%%%%%%%

\setcounter{section}{0}

%%%%%%%%%%%%%%%%%%%%%%%%%%   Section 1 %%%%%%%%%%%%%%%%%%%%%%%%%%%%%%%%%%%%%%%%
\typeout{-----------------------  Section 1  ------------------------}

\tit{Borel Cohomology and Rationalization}

Denote by $\GROUPOIDS$ the category of small groupoids. Let
$\Omega\text{-}\SPECTRA$ be the category of $\Omega$-spectra, where a
morphism $\bff \colon \bfE \to \bfF$ is a sequence of maps
$f_n \colon E_n \to F_n$ compatible with the structure maps and
we work in the category of compactly generated spaces
(see for instance~\cite[Section~1]{Davis-Lueck(1998)}).
A contravariant $\GROUPOIDS$-$\Omega$-spectrum is a contravariant
functor $\bfE \colon \GROUPOIDS \to \Omega\text{-}\SPECTRA$.

Let $\bfE$ be a (non-equivariant) $\Omega$-spectrum. We can associate
to it a contravariant $\GROUPOIDS$-$\Omega$-spectrum
\begin{eqnarray}
\bfE_{\Bor} \colon \GROUPOIDS \to \Omega\text{-}\SPECTRA;
\quad \calg & \mapsto & \map(B\calg;\bfE),
\label{bfE_{Bor}}
\end{eqnarray}
where $B\calg$ is the classifying space associated to $\calg$ and
$\map(B\calg;\bfE)$ is the obvious mapping space spectrum
(see for instance~\cite[page 208 and Definition 3.10 on page
224]{Davis-Lueck(1998)}). In the sequel we use the notion of
an equivariant cohomology theory $\calh^*_?$ with values in $R$-modules
of~\cite[Section 1]{Lueck(2004i)}. It assigns to each (discrete) group $G$
a $G$-cohomology theory $\calh^*_G$ with values in the category of $R$-modules on
the category of pairs of $G$-$CW$-complexes, where $*$ runs through $\bbZ$.
Let $H^*_?(-,\bfE_{\Bor})$ be the to $\bfE_{\Bor}$ associated
equivariant cohomology theory with values in $\bbZ$-modules
satisfying the disjoint union axiom, which has been constructed
in~\cite[Example~1.8]{Lueck(2004i)}. For a given discrete group $G$ and a
$G$-$CW$-pair $(X,A)$ and $n \in \bbZ$ we get a natural identification
\begin{eqnarray}
H^n_G(X,A;\bfE_{\Bor}) & = & H^n(EG \times_G (X,A);\bfE),
\label{H^n_G(X,A;bfE_{Bor})  = H^n(EG times_G (X,A);bfE)}
\end{eqnarray}
where $H^*(-;\bfE)$ is the (non-equivariant) cohomology theory
associated to $\bfE$. It is induced by the following composite of
equivalences of $\Omega$-spectra
\begin{multline*}
\map_{\Or(G)}\left(\map_G(G/?,X)^G,\map\left(B\calg^G(G/H),\bfE\right)\right)
\\
~ \to ~
\map_{\Or(G)}\left(\map_G(G/?,X)^G,\map\left(EG \times_G G/?,\bfE\right)\right)
\\
~ \to ~ \map\left(\map_G(G/?,X)^G \otimes_{\Or(G)} EG \times_G G/?,\bfE\right)
~ \to ~ \map\left(EG \times_G X,\bfE\right)
\end{multline*}
using the notation of~\cite{Lueck(2004i)}. In the literature $H^n(EG \times_G (X,A);\bfE)$
is called \emph{the equivariant Borel cohomology} of $(X,A)$
with respect to the (non-equivariant) cohomology theory $H^*(-;\bfE)$.

Our main example for $\bfE$ will be the topological $K$-theory spectrum $\bfK$,
whose associated (non-equivariant) cohomology theory
$H^*(-;\bfK)$ is topological $K$-theory $K^*$.

There is a functor
$$\Rat\colon \Omega\text{-}\SPECTRA \to  \Omega\text{-}\SPECTRA,
\quad \bfE \mapsto \Rat(\bfE) = \bfE_{(0)},$$
which assigns to an $\Omega$-spectrum $\bfE$ its rationalization $\bfE_{(0)}$.
The homotopy groups $\pi_k(\bfE_{(0)})$ come with a canonical structure of a $\bbQ$-module.
There is a natural  transformation
\begin{eqnarray}
\bfi(\bfE) \colon \bfE & \to & \bfE_{(0)}
\label{bfi(bfE)}
\end{eqnarray}
which induces isomorphisms
\begin{eqnarray}
\pi_k(\bfi(\bfE)) \colon \pi_k(\bfE) \otimes_{\bbZ} \bbQ & \xrightarrow{\cong}&
\pi_k(\bfE_{(0)}).
\label{iso pi_k(bfi(bfE))}
\end{eqnarray}
Composing $\bfE_{\Bor}$ with $\Rat$ yields a contravariant $\Or(G)$-$\Omega$-spectrum
denoted by $\left(\bfE_{\Bor}\right)_{(0)}$. We obtain an equivariant
cohomology theory with values in $\bbQ$-modules by
$H^*_?\left(-;\left(\bfE_{\Bor}\right)_{(0)}\right)$.
The map $\bfi$ induces a natural transformation of equivariant
cohomology theories
\begin{eqnarray}
i^*_?(-;\bfE) \colon H_?^*\left(-;\bfE_{\Bor}\right)
\otimes_{\bbZ} \bbQ
& \to &
H_?^*\left(-;\left(\bfE_{\Bor}\right)_{(0)}\right).
\label{i^*_?(-;bfE)}
\end{eqnarray}

\begin{lemma} \label{lem: i^*_G(X;bfE) is bijective for finite X}
If $G$ is a group $G$ and $(X,A)$ is a relative finite $G$-$CW$-pair, then
$$i^n_G(X,A;\bfE) \colon H_G^n\left(X,A;\bfE_{\Bor}\right) \otimes_{\bbZ} \bbQ \to
H_G^n\left(X,A;\left(\bfE_{\Bor}\right)_{(0)}\right)$$
is a $\bbQ$-isomorphism for all $n \in \bbZ$.
\end{lemma}
\begin{proof}
The transformation $i^*_G(-;\bfE)$ is a natural transformation of $G$-cohomology theories
since $\bbQ$ is flat over $\bbZ$.
One easily checks that it induces a bijection in the case $X = G/H$, since then
there is a commutative square with obvious isomorphisms as vertical
maps and the isomorphism of \eqref{iso pi_k(bfi(bfE))} as lower horizontal arrow
\comsquare{H_G^k\left(G/H;\bfE_{\Bor}\right) \otimes_{\bbZ} \bbQ}{i^*_G(G/H;\bfE)}
{H_G^k\left(G/H;\left(\bfE_{\Bor}\right)_{(0)}\right)}
{\cong}{\cong}
{\pi_{-k}\left(\map(BH,\bfE)\right) \otimes_{\bbZ} \bbQ}{\pi_{-k}\left(\bfi(\map(BH,\bfE))\right)}
{\pi_{-k}\left(\left(\map(BH,\bfE)\right)_{(0)}\right)}
By induction over the number of $G$-cells using Mayer-Vietoris sequences
one shows that $i^*_G(X,A)$ is an isomorphism for all relative finite
$G$-$CW$-pairs $(X,A)$.
\end{proof}

\begin{remark}{\bf (Comparison of the various rationalizations).}
\label{rem: Comparision of the various rationalizations} \em
Notice that $i^*_G(X,A,\bfE)$ of \eqref{i^*_?(-;bfE)} is not an isomorphism for all
$G$-$CW$-pairs $(X,A)$ because the source
does not satisfy the disjoint union axiom for arbitrary index sets in  contrast to the
target. The point is that $- \otimes_{\bbZ} \bbQ$ is compatible with direct sums but not with
direct products.

Since $H_?^*\left(-;\left(\bfK_{\Bor}\right)_{(0)}\right)$ is an
equivariant cohomology theory with values in $\bbQ$-modules satisfying
the disjoint union axiom, we can use the equivariant cohomological
Chern character  of~\cite{Lueck(2004i)} to compute
$H_G^*\left(\underline{E}G;\left(\bfK_{\Bor}\right)_{(0)}\right)$ for all groups
$G$.

This is also true for the
equivariant cohomology theory with values in $\bbQ$-modules satisfying
the disjoint union axiom
$H_?^*\left(-;\left(\bfK_{(0)}\right)_{\Bor}\right)$.
(Here we have changed the order of $\Bor$ and $(0)$.)
But this a much worse approximation of
$K^k(BG) \otimes_{\bbZ} \bbQ$ than
$H_G^*\left(BG;\left(\bfE_{\Bor}\right)_{(0)}\right)$.
Namely, $\bfi$ induces using the universal property of $\Rat$ a natural map of contravariant
$\GROUPOIDS$-$\Omega$-spectra
$$\left(\bfK_{\Bor}\right)_{(0)} \to \left(\bfK_{(0)}\right)_{\Bor}$$
and thus a natural map
$$H_G^k\left(X;\left(\bfK_{\Bor}\right)_{(0)}\right)
\to H_G^k\left(X;\left(\bfK_{(0)}\right)_{\Bor}\right)$$
but this map is in general not an isomorphism. Namely, it is not bijective
for $X = G/H$ for finite non-trivial $H$ and $k = 0$. In this case the source
turns out to be
$$\pi_0\left(\left(\map(BH;\bfK)\right)_{(0)}\right)
~ \cong  ~
K^0(BH) \otimes_{\bbZ} \bbQ
~ \cong ~
\bbQ \times \prod_{p |\; |H|} \left(\bbQ\widehat{_p}\right)^{r(p)}$$
for $r(p)$ the number of conjugacy classes $(h)$  of non-trivial elements $h \in H$ of
$p$-power order, and the target is $K^0(BH;\bbQ)$ which turns out to be isomorphic to $\bbQ$
since the rational cohomology of $BH$ agrees with the one of the one-point-space.
\em
\end{remark}

As mentioned before we want to use the equivariant cohomological Chern character
of~\cite{Lueck(2004i)} to compute
$H_G^*\left(X;\left(\bfK_{\Bor}\right)_{(0)}\right)$.
This requires a careful
analysis of the contravariant functor
$$\FGINJ \to \bbQ\text{ -}\MOD, \quad
H \mapsto H_G^k\left(G/H;\left(\bfK_{\Bor}\right)_{(0)}\right)
= K^k(BH) \otimes_{\bbZ} \bbQ,$$
from the category $\FGINJ$ of finite groups with injective group homomorphisms as
morphisms to the category $\bbQ\text{ -}\MOD$ of $\bbQ$-modules.
It will be carried out in
Section~\ref{sec: The K-Theory of the Classifying Space of a Finite Group}
after some preliminaries in Section~\ref{sec: Some Preliminaries about Pro-Modules}.

%%%%%%%%%%%%%%%%%%%%%%%%%%   Section 2 %%%%%%%%%%%%%%%%%%%%%%%%%%%%%%%%%%%%%%%%

\typeout{-----------------------  Section 2  ------------------------}

\tit{Some Preliminaries about Pro-Modules}

It will be crucial to handle pro-systems and pro-isomorphisms and not
to pass directly to inverse limits.
In this section we fix our notation for handling pro-$R$-modules for a
commutative ring $R$, where ring always means associative ring with unit.
For the definitions in full generality see for instance~\cite[Appendix]{Artin-Mazur(1969)}
or~\cite[\sect~2]{Atiyah-Segal(1969)}.

For simplicity, all pro-$R$-modules dealt
with here will be indexed by the positive  integers.  We
write $\{M_n,\alpha_n\}$ or briefly $\{M_n\}$ for the inverse system
$$ M_0 \xleftarrow{\alpha_1} M_1 \xleftarrow{\alpha_2} M_2 \xleftarrow{\alpha_3}
M_3 \xleftarrow{\alpha_4} \ldots .$$
and also write
$\alpha_n^m := \alpha_{m+1} \circ \cdots \circ \alpha_{n}\colon G_n \to G_m$
for $n > m$ and put $\alpha^n_n =\id_{G_n}$.  For the purposes here, it will
suffice (and greatly simplify the notation)
to work with ``strict'' pro-homomorphisms
$\{f_n\} \colon \{M_n,\alpha_n\} \to \{N_n,\beta_n\}$, i.e.
a collection of  homomorphisms $f_n \colon M_n \to N_n$
for $n \ge 1$ such that  $\beta_{n}\circ f_n = f_{n-1}\circ\alpha_{n}$ holds
for each $ n\ge 2$.  Kernels  and cokernels of strict homomorphisms are
defined in the obvious way.

A pro-$R$-module $\{M_n,\alpha_n\}$  will be called
\emph{pro-trivial} if for each $m \ge 1$, there
is some $n\ge m$ such that  $\alpha^m_n = 0$.  A strict homomorphism
$f\colon \{M_n,\alpha_n\} \to \{N_n,\beta_n\}$ is a
\emph{pro-isomorphism}
if and only if $\ker(f)$ and $\cok(f)$ are both
pro-trivial, or, equivalently, for each $m\ge 1$ there is some
$n\ge m$ such that $\im(\beta_n^m) \subseteq \im(f_m)$
and $\ker(f_n) \subseteq \ker(\alpha_n^m)$.
A sequence of strict homomorphisms
$$\{M_n,\alpha_n\} \xrightarrow{\{f_n\}} \{M_n',\alpha_n'\}
\xrightarrow{g_n} \{M_n'',\alpha_n''\}$$
will be called \emph{exact} if the sequences of $R$-modules
$M_n \xrightarrow{f_n} N_n \xrightarrow{g_n} M_n''$ is
exact for each $n \ge 1$, and it is called \emph{pro-exact} if
$g_n \circ f_n = 0$  holds for $n \ge 1$ and
the pro-$R$-module  $\{\ker(g_n)/\im(f_n)\bigr\}$ is pro-trivial.

The following results will be needed later.

\begin{lemma} \label{lem: pro-exactness and  limits}
Let $0 \to \{M_n',\alpha_n'\} \xrightarrow{\{f_n\}} \{M_n,\alpha_n\}
\xrightarrow{\{g_n\}} \{M_n'',\alpha_n''\} \to 0$ be a pro-exact sequence
of
pro-$R$-modules. Then there is a natural exact sequence
\begin{multline*}
0 \to \invlim{n \ge 1}{M_n'} \xrightarrow{\invlim{n \ge 1}{f_n}}
\invlim{n \ge 1}{M_n} \xrightarrow{\invlim{n \ge 1}{g_n}}
\invlim{n \ge 1}{M_n''} \xrightarrow{\delta}
\\
 \higherlim{n \ge 1}{1}{M_n'} \xrightarrow{\higherlim{n \ge 1}{1}{f_n}}
\higherlim{n \ge 1}{1}{M_n} \xrightarrow{\higherlim{n \ge 1}{1}{g_n}}
\higherlim{n \ge 1}{1}{M_n''} \to 0.
\end{multline*}
In particular a pro-isomorphism
$\{f_n\} \colon \{M_n,\alpha_n\} \to \{N_n,\beta_n\}$ induces isomorphisms
$$\begin{array}{llcl}
\invlim{n \ge 1}{f_n} \colon & \invlim{n \ge 1}{M_n}
& \xrightarrow{\cong} & \invlim{n \ge 1}{N_n};
\\
\higherlim{n \ge 1}{1}{f_n} \colon & \higherlim{n \ge 1}{1}{M_n}
& \xrightarrow{\cong} & \higherlim{n \ge 1}{1}{N_n}.
\end{array}$$
\end{lemma}
\begin{proof} If $0 \to \{M_n',\alpha_n'\} \xrightarrow{\{f_n\}} \{M_n,\alpha_n\}
\xrightarrow{g_n} \{M_n'',\alpha_n''\} \to 0$ is exact, the
construction of the six-term sequence is standard
(see for instance~\cite[Proposition~7.63 on page~127]{Switzer(1975)}).
Hence it remains to show for a pro-trivial pro-$R$-module
$\{M_n,\alpha_n\}$ that $\invlim{n \ge 1}{M_n}$ and
$\higherlim{n \ge 1}{1}{M_n}$ vanish. This follows directly from the
standard construction for these limits as the kernel and cokernel of
the homomorphism
$$\prod_{n \ge 1} M_n \to \prod_{n \ge 1} M_n, \quad
(x_n)_{n \ge 1} ~ \mapsto
(x_n - \alpha_{n+1}(x_{n+1}))_{n \ge 1}. $$
\end{proof}

\begin{lemma} \label{lem: pro-exactness and exactness}
Fix any commutative Noetherian ring $R$, and any ideal
$I\subseteq R $.  Then for any exact sequence $M' \to M \to M''$ of
finitely generated $R$-modules, the sequence
        $$ \{M'/I^nM'\} \to  \{M/I^nM\} \to  \{M''/I^nM''\} $$
of pro-$R$-modules is pro-exact.
\end{lemma}
\begin{proof}  It suffices to prove this for a short exact sequence
$0 \to M'\to M\to M'' \to 0$. Regard
$M'$ as a submodule of $M$, and consider the exact sequence
$$ 0 \to \left\{{\tfrac{(I^n M){\cap} M'}{I^n M'}}\right\}
\to  \{M'/I^nM'\} \to  \{M/I^nM\} \to  \{M''/I^nM''\}
        \to 0. $$
By~\cite[Theorem~10.11 on page~107]{Atiyah-McDonald(1969)},
the filtrations $\{(I^nM){\cap}M'\}$ and
$\{I^nM'\}$ of $M'$ have ``bounded difference'', i.e. there exists
$k>0$ with the property that
$(I^{n+k}M){\cap}M'\subseteq I^nM'$
holds for all $n \ge 1$.  The first
term in the above exact sequence is thus pro-trivial, and so the remaining
terms define a short sequence of pro-$R$-modules which is pro-exact.
\end{proof}

%%%%%%%%%%%%%%%%%%%%%%%%%%   Section 3 %%%%%%%%%%%%%%%%%%%%%%%%%%%%%%%%%%%%%%%%

\typeout{-----------------------  Section 3  ------------------------}

\tit{The K-Theory of the Classifying Space of a Finite Group}

Next we investigate the contravariant functor
from the category $\FGINJ$ of finite groups with injective group homomorphisms as
morphisms to the category $\bbZ\text{ -}\MOD$ of $\bbZ$-modules
$$\FGINJ \to \bbZ\text{ -}\MOD, \quad
H \mapsto K^k(BH).$$

We need some input from representation theory.
Recall that $R(G)$ denotes the complex representation ring.
Let $\bbI_G$ be the kernel of the ring homomorphism
$\res_G^{\{1\}}\colon R(G) \to R(\{1\}$ which is
the same as the kernel of augmentation
ring homomorphism $R(G) \to \bbZ$ sending $[V]$ to $\dim_{\bbC}(V)$.  We will frequently
use the so called \emph{double coset formula} (see
\cite[Proposition~22 in Chapter~7 on page~58]{Serre(1977)}).
It says for two subgroups $H,K \subseteq G$
\begin{eqnarray}
\res_G^K \circ \ind_H^G & =  & \sum_{KgH \in K\backslash G/H}
\ind_{c(g)\colon H\cap g^{-1}Kg \to K}
 \circ \res_{H}^{ H\cap g^{-1}Kg},
\label{double coset formula}
\end{eqnarray}
where $c(g)$ is conjugation with $g$, i.e. $c(g)(h) = ghg^{-1}$, and
$\ind$ and $\res$ denote induction and restriction. One consequence
of it is that $\ind_H^G\colon  R(H) \to R(G)$ sends $\bbI_H$ to $\bbI_G$.
Obviously $\res_G^H\colon R(G) \to R(H)$ maps
$\bbI_G$ to $\bbI_H$.

For an abelian group $M$ let $M_{(p)}$ be the localization of $M$ at $p$.
If $\bbZ_{(p)}$ is the subring of $\bbQ$ obtained from $\bbZ$
by inverting all prime numbers except $p$,
then $M_{(p)} = M \otimes_{\bbZ} \bbZ_{(p)}$. Recall that the functor $ ? \otimes_{\bbZ}
\bbZ_{(p)}$ is exact.

\begin{lemma} \label{lem: res circ ind and res have the same image}
Let $G$ be a finite group. Let $p$ be a prime number and denote by
$G_p$ a $p$-Sylow subgroup of $G$. Then the composite
$$R(G_p)_{(p)} \xrightarrow{\ind_{G_p}^G} R(G)_{(p)}
\xrightarrow{\res_G^{G_p}} R(G_p)_{(p)}$$
has the same image as
$$\res_G^{G_p} \colon R(G)_{(p)} \to R(G_p)_{(p)}.$$
\end{lemma}
\begin{proof} A subgroup $H \subseteq G$ is called \emph{$p$-elementary} if it is
isomorphic to $C \times P$ for a cyclic group $C$
of order prime to $p$ and a $p$-group $P$.
Let $\{C_i \times P_i \mid i = 1,2, \ldots , r\}$
be a complete system of representatives
of conjugacy classes of $p$-elementary subgroups of $G$. We can assume without loss of generality $P_i \subseteq G_p$.
Define for $i = 1,2,  \ldots , r$
a homomorphism of abelian groups
$$\phi_i ~ := ~
\sum_{\substack{G_p \cdot g \cdot (C_i \times P_i)
\in \\ G_p\backslash G/(C_i \times P_i)}}
\ind_{c(g)\colon P_i \cap g^{-1}G_pg \to G_p} \circ \res_{P_i}^{P_i \cap g^{-1}G_pg}
\colon R(P_i) ~ \to ~ R(G_p).$$
Since the order of $C_i$ is prime to $p$, we have
$(C_i \times P_i) \cap g^{-1}G_p g = P_i \cap g^{-1}G_p g$ for $g \in G$. Hence the
following diagram commutes (actually before localization)
by the double coset formula
$$\begin{CD}
\bigoplus_{i=1}^r R(P_i)_{(p)}
@> \bigoplus_{i=1}^r \ind_{P_i}^{G_p} >>
R(G_p)_{(p)}
\\
@V \bigoplus_{i=1}^r \ind_{P_i}^{C_i \times P_i} VV         @VV \ind_{G_p}^G V
\\
\bigoplus_{i=1}^r R(C_i \times P_i)_{(p)}
@> \bigoplus_{i=1}^r \ind_{C_i \times P_i}^{G_p} >>
R(G)_{(p)}
\\
@V \bigoplus_{i=1}^r \res_{C_i \times P_i}^{P_i} VV         @VV \res_G^{G_p} V
\\
\bigoplus_{i=1}^r R(P_i)_{(p)}
@> \bigoplus_{i=1}^r \phi_i >>
R(G_p)_{(p)}
\end{CD}$$
The middle horizontal arrow $\bigoplus_{i=1}^r \ind_{C_i \times P_i}^{G_p}$
is surjective by Brauer's Theorem
\cite[Theorem~18 in Chapter~10 on page~75]{Serre(1977)}. The composite
of the left lower vertical arrow and the left upper vertical arrow
$\bigoplus_{i=1}^r \res_{C_i \times P_i}^{P_i} \circ \ind_{P_i}^{C_i \times P_i}$
is $\bigoplus_{i=1}^r |C_i| \cdot \id$ and hence an isomorphism.
Now the claim follows from an easy diagram chase.
\end{proof}

\begin{lemma} \label{lem: res circ ind for different primes}
Let $p$ and $q$ be different primes. Then the composition
$$R(G_p) \xrightarrow{\ind_{G_p}^G} R(G)
\xrightarrow{\res_G^{G_q}} R(G_q)$$
agrees with $|G_q\backslash G / G_p| \cdot
\ind_{\{1\}}^{G_q} \circ \res_{G_p}^{\{1\}}.$
\end{lemma}
\begin{proof} This follows from the double coset formula \eqref{double coset formula} since
$G_p \cap g^{-1}G_qg = \{1\}$ for each  $g \in G$.
\end{proof}

\begin{lemma} \label{lem: computing R(G) to prod_p R(G_p)}
Let $G$ be a finite group and let $\bbI_G \subseteq R(G)$
be the augmentation ideal. Then
the following sequence of $R(G)$-modules is exact
$$0 \to \bigcap_{m \ge 1} (\bbI_G)^m \xrightarrow{i} \bbI_G
\xrightarrow{\prod_{p} \res_{G}^{G_p}} \prod_{p \in \calp(G)} \im\left(\res_{G}^{G_p}\colon \bbI_G
  \to \bbI_{G_p}\right) \to 0,$$
where $i$ is the inclusion and $\calp(G)$ is the set of primes dividing $|G|$.
\end{lemma}
\begin{proof} The kernel of $\prod_{p} \res_G^{G_p} \colon R(G) \to \prod_{p} R(G_p)$
is $\bigcap_{m \ge 1} (\bbI_G)^m$  by~\cite[Proposition~6.12 on page~269]{Atiyah(1961)}.
Hence it remains to show that
$$\prod_{p} \res_{G}^{G_p} \colon \bbI_G
\to \prod_{p} \im\left(\res_{G}^{G_p}\colon \bbI_G \to \bbI_{G_p}\right)$$
is surjective. It suffices to show for a each  prime number $q$ that its localization
$$\prod_{p} \res_{G}^{G_p} \colon (\bbI_G)_{(q)}
\to \prod_{p} \im\left(\res_{G}^{G_p}\colon (\bbI_G)_{(q)}  \to (\bbI_{G_p})_{(q)} \right)$$
is surjective. Next we construct the following commutative diagram
$$\begin{CD}
\bigoplus_{p \not= q} (\bbI_{G_p})_{(q)}
@> \prod_{p \not= q} \res_G^{G_p} \circ \ind_{G_p}^G >>
\prod_{p \not=q} \im\left(\res_G^{G_p} \colon (\bbI_G)_{(q)} \to (\bbI_{G_p})_{(q)} \right)
\\
@V \bigoplus_{p \not= q} \ind_{G_p}^{G} VV @V i VV
\\
(\bbI_G)_{(q)}
@> \prod_{p} \res_G^{G_p} >>
\prod_{p} \im\left(\res_G^{G_p} \colon (\bbI_G)_{(q)} \to (\bbI_{G_p})_{(q)} \right)
\\
@V p_1 VV  @V p_2 VV
\\
\cok\left(\bigoplus_{p \not= q} \ind_{G_p}^{G}\right)
@> f >>
\im\left(\res_G^{G_q} \colon (\bbI_G)_{(q)} \to (\bbI_{G_q})_{(q)} \right)
\end{CD}$$
Here $i$ is the inclusion and $p_1$ and $p_2$ are the obvious projections.
Since the composition
\begin{multline*}
\bigoplus_{p \not= q} (\bbI_{G_p})_{(q)}
\xrightarrow{ \bigoplus_{p \not= q} \ind_{G_p}^{G}}
(\bbI_G)_{(q)}
\xrightarrow{\prod_{p} \res_G^{G_p}}
\prod_{p} \im\left(\res_G^{G_p} \colon (\bbI_G)_{(q)} \to (\bbI_{G_p})_{(q)} \right)
\\
\xrightarrow{p_2}
\im\left(\res_G^{G_q} \colon (\bbI_G)_{(q)} \to (\bbI_{G_q})_{(q)} \right)
\end{multline*}
agrees with
$$\bigoplus_{p \not= q} \res_{G}^{G_q} \circ \ind_{G_p}^G \colon
\bigoplus_{p \not= q} (\bbI_{G_p})_{(q)}
\to \im\left(\res_G^{G_q} \colon (\bbI_G)_{(q)} \to (\bbI_{G_q})_{(q)} \right)$$
and hence is trivial by Lemma~\ref{lem: res circ ind for different primes},
there exists a map
$$f\colon \cok\left(\bigoplus_{p \not= q} \ind_{G_p}^{G}\right) \to
\im\left(\res_G^{G_q} \colon (\bbI_G)_{(q)} \to (\bbI_{G_q})_{(q)} \right)$$
such that the diagram above commutes. Since
$$p_2 \circ \prod_{p} \res_G^{G_p} = \res_G^{G_q}\colon (\bbI_G)_{(q)} \to
\im\left(\res_G^{G_q} \colon (\bbI_G)_{(q)} \to (\bbI_{G_q})_{(q)} \right)$$
is by definition surjective, $f$ is surjective. The upper horizontal arrow in the
commutative diagram above is surjective by Lemma
\ref{lem: res circ ind and res have the same image}. Now the claim follows by an easy
diagram chase.
\end{proof}

\begin{theorem}[Structure of $\{\bbI_G/(\bbI_G)^{n+1}\}$]
\label{the: computation of AI_G/(AI_G)^{n+1}}
Let $G$ be a finite group. Let $\calp(G)$ be the set of primes dividing $|G|$.

\begin{enumerate}

\item \label{the: computation of AI_G/(AI_G)^{n+1}: p^a I and i^2}
There are positive integers $a$, $b$  and $c$
such that for each prime $p$ dividing the order of $|G|$
\begin{eqnarray*}
p^a \cdot \bbI_{G_p} & \subseteq & \bbI_{G_p}^2;
\\
\bbI_{G_p}^b & \subseteq & p \cdot \bbI_{G_p};
\\
\bbI_G \cdot \bbI_{G_p} & \subseteq & \bbI_{G_p}^2;
\\
(\bbI_{G_p})^c  & \subseteq & \bbI_G \cdot \bbI_{G_p};
\end{eqnarray*}

\item \label{the: computation of AI_G/(AI_G)^{n+1}: computing I/I^n}
For a prime $p$ dividing $|G|$ let $\im(\res_G^{G_p})$ be the image of
$\res_G^{G_p}\colon \bbI_G \to \bbI_{G_p}$.
We obtain a sequence of pro-isomorphisms of pro-$\bbZ$-modules
\begin{multline*}
\{\bbI_G/(\bbI_G)^{n+1}\}
\xrightarrow{\cong}
\prod_{p \in \calp(G)} \{\im(\res_G^{G_p})/(\bbI_G)^n \cdot \im(\res_G^{G_p})\}
\\
\xrightarrow{\cong}
\prod_{p \in \calp(G)} \{\im(\res_G^{G_p})/(\bbI_{G_p})^n \cdot \im(\res_G^{G_p})\}
\\
\xleftarrow{\cong}
\prod_{p \in \calp(G)} \{\im(\res_G^{G_p})/(\bbI_{G_p})^{bn}\cdot \im(\res_G^{G_p})\}
\\
\xrightarrow{\cong}
\prod_{p \in \calp(G)} \{\im(\res_G^{G_p})/p^n \cdot \im(\res_G^{G_p})\} .
\end{multline*}

\item \label{the: computation of AI_G/(AI_G)^{n+1}: computing R/I^n}
There is an isomorphism of pro-$\bbZ$-modules
$$\{\bbZ\} \oplus \{\bbI_G/(\bbI_G)^{n}\}
\xrightarrow{\cong} \{R(G)/(\bbI_G)^n\},$$
where $\{\bbZ\}$ denotes the constant inverse system
 $\bbZ \xleftarrow{\id} \bbZ \xleftarrow{\id} \ldots.$
\end{enumerate}
\end{theorem}
\begin{proof}
\eqref{the: computation of AI_G/(AI_G)^{n+1}: p^a I and i^2}
The existence of the integers $a$, $b$ and
$c$ for which the inclusions appearing in the
statement of Theorem~\ref{the: computation of AI_G/(AI_G)^{n+1}} hold follows from
results of~\cite[Theorem~6.1 on page~265]{Atiyah(1961)} and
\cite[Proposition~1.1 in Part~III on page~277]{Atiyah-Tall(1969)}.
\\[1mm]
\eqref{the: computation of AI_G/(AI_G)^{n+1}: computing I/I^n}
These inequalities of assertion
\eqref{the: computation of AI_G/(AI_G)^{n+1}: p^a I and i^2} imply
that the second, third and fourth map of pro-$\bbZ$-isomorphism appearing in the
statement of Theorem~\ref{the: computation of AI_G/(AI_G)^{n+1}} are indeed well-defined
pro-isomorphisms. The first map
$$\{\bbI_G/(\bbI_G)^{n+1}\} \xrightarrow{\cong}
\prod_{p \in \calp(G)} \{\im(\res_G^{G_p})/(\bbI_G)^n \cdot \im(\res_G^{G_p})\}$$
is a well-defined pro-isomorphism of pro-$\bbZ$-modules by
Lemma~\ref{lem: pro-exactness and exactness} and
Lemma~\ref{lem: computing R(G) to prod_p R(G_p)} provided
$\left\{\left(\bigcap_{m \ge 1} (\bbI_G)^m\right)/\bbI_G^n \cdot  \left(\bigcap_{m \ge 1} (\bbI_G)^m\right)\right\}$ is pro-trivial.
The latter statement follows from Lemma~\ref{lem: pro-exactness and exactness}
applied to the exact sequence
\begin{multline*}
0 \to \bigcap_{m \ge 1} (\bbI_G)^m \to \bbI_G \to \bbI_G/\bigcap_{m \ge 1} (\bbI_G)^m \to 0.
\end{multline*}
\\[1mm]
\eqref{the: computation of AI_G/(AI_G)^{n+1}: computing R/I^n}
Consider the isomorphism of finitely generated free abelian groups
$$ \bbZ \oplus \bbI_G \xrightarrow{\cong} R(G),
\quad (m,x) \mapsto x + m \cdot [\bbC].$$
It becomes an isomorphism of rings if we equip the source with the multiplication
$(m,x) \cdot (n,y) = (mn,my + nx + xy)$. In particular
$\bbI_G^n  \cdot  ( \bbI_G \oplus\bbZ ) \subseteq \bbI_G^n\oplus 0$ for $n \ge 1$.
This finishes the proof of Theorem~\ref{the: computation of AI_G/(AI_G)^{n+1}}.

\end{proof}

Now we can give the proof of Theorem~\ref{the: computation of K_*(BG) and K^*(BG) for finite G}.
\begin{proof} In the sequel we abbreviate
$\im(\res_G^{G_p}) =  \im\left(\res_G^{G_p}\colon \bbI_G
  \to \bbI_{G_p}\right)$. Notice that
$\im(\res_G^{G_p}) \subseteq R(G_p)$ is a finitely
generated free $\bbZ$-module. We obtain
from Lemma~\ref{lem: pro-exactness and  limits}  and
Theorem~\ref{the: computation of AI_G/(AI_G)^{n+1}} an isomorphism
$$
\invlim{n \ge 1}{R(G)/(\bbI_G)^n} ~ \cong ~
\bbZ \times \prod_{p \in \calp(G)} \invlim{n \ge 1}{\im(\res_G^{G_p})/p^n\cdot \im(\res_G^{G_p})}.$$
Now the Atiyah-Segal Completion Theorem ~\cite{Atiyah-Segal(1969)}
yields an isomorphisms
\begin{eqnarray}
\invlim{n \ge 1}{R(G)/(\bbI_G)^n} ~ \xrightarrow{\cong}  ~
\invlim{n \ge 1} K^0((BG)_n) ~ \xleftarrow{\cong} ~ K^0(BG)
& & \label{AS for finite groups}
\end{eqnarray}
and $K^1(BG) = 0$. This implies
\begin{eqnarray*}
K^0(BG) & \cong & \bbZ \oplus \bigoplus_{p \in \calp(G)} \im\left(\res_G^{G_p}\colon \bbI_G
  \to \bbI_{G_p}\right) \otimes_{\bbZ} \bbZ\widehat{_p};
\\
K^1(BG) & \cong & 0.
\end{eqnarray*}

Next we show that the rank of the finitely generated free abelian group
$\im\left(\res_G^{G_p}\colon \bbI_G   \to \bbI_{G_p}\right) \subseteq R(G_p)$
is the number $r(p)$ of conjugacy classes $(g)$
of elements $g \in G$ whose order $|g|$ is $p^d$ for some integer $d \ge 1$.
This follows from the commutative diagram
\comsquare{\bbC \otimes_{\bbZ} R(G)}{\res_G^{G_p}}{\bbC \otimes_{\bbZ} R(G_p)}
{\cong}{\cong}{\class_{\bbC}(G)}{\res_G^{G_p}}{\class_{\bbC}(G_p)}
where $\class_{\bbC}(G)$ denotes the complex vector space of class functions on $G$, i.e.
functions $G \to \bbC$ which are constant on
conjugacy classes of elements, (and analogous for $G_p$),
the vertical isomorphisms are given by taking the character of a complex representation,
and the lower horizontal arrow is given by restricting a function $G \to \bbC$ to $G_p$.

Recall that $\bbI_p(G)$ is canonically isomorphic to
$\im\left(\res_G^{G_p}\colon \bbI_G  \to \bbI_{G_p}\right)$.

One easily checks that the isomorphisms obtained from the one appearing
in Theorem~\ref{the: computation of AI_G/(AI_G)^{n+1}}
\eqref{the: computation of AI_G/(AI_G)^{n+1}: computing I/I^n} and
\eqref{the: computation of AI_G/(AI_G)^{n+1}: computing R/I^n}
by applying the inverse limit and the isomorphism
\eqref{AS for finite groups} are compatible with the obvious multiplicative structures.

This finishes the proof of
Theorem~\ref{the: computation of K_*(BG) and K^*(BG) for finite G}.
\end{proof}

%%%%%%%%%%%%%%%%%%%%%%%%%%% Section 4 %%%%%%%%%%%%%%%%%%%%%%%%%%%%%%%%%%
\typeout{-----------------------  Section 4  ------------------------}

\tit{Proof of the Main Result}

In this section we want to prove our main Theorem~\ref{the: main theorem}.
We want to apply the cohomological equivariant Chern character
of~\cite{Lueck(2004i)} to the equivariant
cohomology theory $H^*_?\left(-;\left(\bfK_{\Bor}\right)_{(0)}\right)$.
This requires to analyze the contravariant functor
\begin{eqnarray}
\FGINJ & \to & \bbQ\text{ -}\MOD, \quad
H \mapsto H_G^l\left(G/H;\left(\bfK_{\Bor}\right)_{(0)}\right).
\label{functor H_G^l(G/H;(bfK_{Bor})_{(0)})}
\end{eqnarray}
From \eqref{H^n_G(X,A;bfE_{Bor})  = H^n(EG times_G (X,A);bfE)}
and Lemma~\ref{lem: i^*_G(X;bfE) is bijective for finite X}
we conclude that  the contravariant functor \eqref{functor H_G^l(G/H;(bfK_{Bor})_{(0)})}
is naturally equivalent to the contravariant functor
\begin{eqnarray}
\FGINJ & \to & \bbQ\text{ -}\MOD, \quad
H \mapsto K^l(BH) \otimes_{\bbZ} \bbQ.
\label{functor  K^l(BH) otimes_Z Q}
\end{eqnarray}
Theorem~\ref{the: computation of K_*(BG) and K^*(BG) for finite G} yields
the contravariant functor \eqref{functor  K^l(BH) otimes_Z Q}
is trivial for odd $l$ and is naturally equivalent to the contravariant functor
\begin{eqnarray}
\FGINJ &\to &\bbQ\text{ -}\MOD
\quad H \mapsto \bbQ \times \prod_{p} \bbI_p(H) \otimes_{\bbZ} \bbQ\widehat{_p}
\label{splitting of K^l(B?) otimes_{bbZ}  bbQ}
\end{eqnarray}
 for even $l$,  where the factor $\bbQ$ is constant in $H$ and functoriality for the other factors
is given by restriction.

Given a contravariant functor $F \colon \FGINJ \to \bbQ\text{ -}\MOD$, define
the $\bbQ[\aut(H)]$-module
\begin{eqnarray}
T_H F(H) & := &
\ker\left(\prod_{K \subsetneq H} F(K \hookrightarrow H) \colon
F(H) ~ \to ~ \prod_{K \subsetneq H} F(K) \right).
\label{def of T_HF(H)}
\end{eqnarray}
Next we compute $T_H\left(K^0(BH) \otimes_{\bbZ} \bbQ\right)$.
Since $T_H$ is compatible with direct products, we obtain from
\eqref{splitting of K^l(B?) otimes_{bbZ}  bbQ} a canonical $\bbQ[\aut(H)]$-isomorphism
\begin{eqnarray}
T_H \left(K^0(BH) \otimes_{\bbZ} \bbQ\right) & = &
T_H(\bbQ) \times \prod_{p} T_H\left(\bbI_p(H) \otimes_{\bbZ} \bbQ\widehat{_p}\right).
\label{splitting of T_HK^l(BH)}
\end{eqnarray}
Since $\bbQ$ is the constant functor, we get
\begin{eqnarray}
T_H(\bbQ) ~ := ~ \left\{\begin{array}{lll}
0& & \text{ if } H \not= \{1\};
\\
 \bbQ & & \text{ if } H = \{1\}.
\end{array}\right. \label{computation of T_H(qq)}
\end{eqnarray}
Fix a prime number $p$.
Since for any finite group $H$ the map given by restriction to finite cyclic subgroups
$$R(H) \to \prod_{\substack{C \subseteq H\\C \text{ cyclic}}} R(C)$$
is injective, we conclude
\begin {lemma} \label{lem: T_H(I_p(H) otimes_Z Q_p) for H non finite cyclic p-group}
For a finite group $H$
$$T_H\left(\bbI_p(H)\right) ~ = ~ 0,$$
unless $H$ is a non-trivial cyclic $p$-group.
\end{lemma}

Let $C$ be a non-trivial finite cyclic
$p$-group. Then we get
\begin{eqnarray}
T_C \left(\bbI_p(C)\right) & = &
\ker\left(\res_C^{C'} \colon R(C) \to R(C')\right),
\label{computation of T_C (bbI_p(C)) for C Z/p^k}
\end{eqnarray}
where $C' \subseteq C$ is the unique cyclic subgroup of index $p$ in $C$.

Recall that taking the character of a rational
representation of a finite group $H$ yields an isomorphism
$$\chi \colon R_{\bbQ}(H) \otimes_{\bbZ} \bbQ \xrightarrow{\cong} \class_{\bbQ}(H),$$
where $R_{\bbQ}(C)$ is the rational representation ring of $C$ and
$\class_{\bbQ}(H)$ is the rational vector space of functions $f\colon H \to \bbQ$
for which $f(g_1) = f(g_2)$ holds if the cyclic subgroups generated by $g_1$ and $g_2$
are conjugate in $H$ (see~\cite[page~68 and Theorem~29 on page~102]{Serre(1977)}).
Hence there is an idempotent
$\theta_C \in  R_{\bbQ}(C)  \otimes_{\bbZ} \bbQ$ which is uniquely determined by the property
that its character sends a generator of $C$ to $1$ and all other elements to $0$.
Denote its image under the change of coefficients map
$R_{\bbQ}(C)  \otimes_{\bbZ} \bbQ  \to R(C)  \otimes_{\bbZ} \bbQ$ also by $\theta_C$.
Let $ \theta_C \cdot R(C) \otimes_{\bbZ} \bbQ\widehat{_p} \subseteq R(C) \otimes_{\bbZ} \bbQ\widehat{_p}$
be the image of the idempotent endomorphism
$R(C) \otimes_{\bbZ} \bbQ\widehat{_p} \to R(C) \otimes_{\bbZ} \bbQ\widehat{_p}$ given by multiplication with
$\theta_C$.

\begin{lemma}
 \label{lem: T_C(I_p(C) otimes_Z Q_p) for C a finite cyclic p-group}
For every non-trivial cyclic $p$-group $C$ the inclusion induces a $\bbQ[\aut(C)]$-isomorphism
$$\theta_C \cdot R(C) \otimes_{\bbZ} \bbQ ~ \xrightarrow{\cong} ~
T_C \left(\bbI_p(C) \otimes_{\bbZ} \bbQ\right).$$
\end{lemma}
\begin{proof}
Since the map
$\res_C^{C'} \colon R(C) \otimes_{\bbZ} \bbQ \to R(C') \otimes_{\bbZ} \bbQ$
sends $\theta_C$ to zero, $\theta_C \cdot R(C) \otimes_{\bbZ} \bbQ$ is contained in
$\ker\left(\res_C^{C'} \colon R(C) \to R(C')\right) \otimes_{\bbZ} \bbQ$.
For $x \in \ker\left(\res_C^{C'} \colon R(C) \to R(C')\right) \otimes_{\bbZ} \bbQ$
one gets $\theta_C \cdot x - x = 0$ by the  calculation appearing
in the proof of~\cite[Lemma~3.4~(b)]{Lueck(2002d)}.
\end{proof}

\begin{lemma} \label{lem: rational computation of H^*_G(X;(bfK_{Bor})_{(0)})}
For every proper $G$-$CW$-complex $X$ and $n \in \bbZ$ there is an isomorphism,
natural in $X$,
\begin{multline*}
\overline{\ch}_G^n \colon H^*_G\left(X;\left(\bfK_{\Bor}\right)_{(0)}\right)    \xrightarrow{\cong}
\\
\prod_{i \in \bbZ} H^{2i+n}(G\backslash X;\bbQ) ~ \times ~
\prod_{p} ~ \prod_{(C)\in \calc_p(G)} ~
\prod_{i \in \bbZ} H^{2i+n}_ {W_GC}(C_GC \backslash X^C;\theta_C \cdot R(C) \otimes_{\bbZ} \bbQ\widehat{_p}),
\end{multline*}
where $\calc_p(C)$ is the set of conjugacy classes of non-trivial cyclic $p$-subgroups of $G$
and $W_GC = N_GC/C_GC$ is considered as a subgroup of $\aut(C)$ and thus acts on
$\theta_C \cdot R(C) \otimes_{\bbZ} \bbQ\widehat{_p}$.
\end{lemma}
\begin{proof}
This follows from~\cite[Theorem~5.5~(c) and Example~5.6]{Lueck(2004i)} using
\eqref{computation of T_H(qq)}, Lemma~\ref{lem: T_H(I_p(H) otimes_Z Q_p) for H non finite cyclic p-group}
and Lemma~\ref{lem: T_C(I_p(C) otimes_Z Q_p) for C a finite cyclic p-group}.
\end{proof}

For a generator $t \in C$ let $\bbC_t$ be the $\bbC$-representation with
$\bbC$ as underlying complex vector space such that $t$ operates on $\bbC$ by multiplication
with $\exp\left(\frac{2\pi i}{|C|}\right)$. Let $\Gen(C)$ be the set of generators.
Notice that $\aut(C)$ acts in an obvious way on $\Gen(C)$ such that the $\aut(C)$-action is transitive and free,
and acts on $R(C)$ by restriction. In the sequel $\chi_V$ denotes for
a complex representation $V$ its character.

\begin{lemma} \label{lem: identifying Q[Gen(C)] and theta_C cdot R(C) otimes_Z Q)}
Let $C$ be a finite cyclic group. Then

\begin{enumerate}

\item \label{lem: identifying Q[Gen(C)] and theta_C cdot R(C) otimes_Z Q): C}
The map
$$v(C) \colon \theta_C \cdot R(C) \otimes_{\bbZ} \bbC ~ \xrightarrow{\cong} ~
\prod_{\Gen(C)} \bbC,
\quad [V] \mapsto \left(\chi_V(t)\right)_{t \in \Gen(C)}
$$
is a $\bbC[\aut(C)]$-isomorphism if $\aut(C)$ acts on the target by permuting the
factors.
The map $v(C)$ is compatible with the ring structure
on the source induced by the tensor product of representations and the
product ring structure on the target;

\item \label{lem: identifying Q[Gen(C)] and theta_C cdot R(C) otimes_Z Q): Q}
There is an isomorphism of $\bbQ[\aut(C)]$-modules
$$u(C) \colon \bbQ[\Gen(C)] ~ \xrightarrow{\cong} ~
\theta_C \cdot R(C) \otimes_{\bbZ} \bbQ.$$
\end{enumerate}
\end{lemma}
\begin{proof}
\eqref{lem: identifying Q[Gen(C)] and theta_C cdot R(C) otimes_Z Q): C}
The map
$$R(C) \otimes_{\bbZ} \bbC \xrightarrow{\cong} \prod_{g \in C} \bbC, \quad [V] \mapsto (\chi_V(g))_{g \in C}$$
is an isomorphism of rings. One easily checks that it is compatible with the
$\aut(C)$ actions. Now the assertion follows from the fact that the character of $\theta_C$
sends a generator of $C$ to $1$ and any other element of $C$ to $0$.
\\[1mm]
\eqref{lem: identifying Q[Gen(C)] and theta_C cdot R(C) otimes_Z Q): Q}
Obviously $\bbQ[\Gen(C)]$ is $\bbQ[\aut(C)]$-isomorphic to the regular
representation $\bbQ[\aut(C)]$ since $\Gen(C)$ is a transitive free $\aut(C)$-set.
It remains to show that $\theta_C \cdot R(C) \otimes_{\bbZ} \bbQ$
is $\bbQ[\aut(C)]$-isomorphic to the regular
representation $\bbQ[\aut(C)]$. By character theory it suffices to show that
$\theta_C \cdot R(C) \otimes_{\bbZ} \bbC$
is $\bbC[\aut(C)]$-isomorphic to the regular
representation $\bbC[\aut(C)]$. This follows from assertion
\eqref{lem: identifying Q[Gen(C)] and theta_C cdot R(C) otimes_Z Q): C}.
\end{proof}

\begin{lemma} \label{lem: refined non-multiplicative rational computation of H^*_G(X;(bfK_{Bor})_{(0)})}
For every proper $G$-$CW$-complex $X$ and $n \in \bbZ$ there is an isomorphism,
natural in $X$,
\begin{multline*}
\overline{\overline{\ch}}_G^n \colon H^*_G\left(X;\left(\bfK_{\Bor}\right)_{(0)}\right)    \xrightarrow{\cong}
\\
\prod_{i \in \bbZ} H^{2i+n}(G\backslash X;\bbQ) \times
\prod_{p} \prod_{(g) \in \con_p(G)}
H^{2i+n}(C_G\langle g \rangle\backslash X^{\langle g \rangle};\bbQ\widehat{_p}).
\end{multline*}
\end{lemma}
\begin{proof} Fix a prime $p$.
Let $C$ be a cyclic subgroup of $G$ of order $p^d$ for some integer $d \ge 1$.
The obvious $N_GC$-action on $C$ given by conjugation induces an embedding of groups
$W_GC \to \aut(C)$. The obvious action of $\aut(C)$ on $\Gen(C)$ is free and transitive. Thus we obtain an isomorphism
of $\bbQ\widehat{_p}[W_GC]$-modules
$$\bbQ\widehat{_p}[\Gen(C)] \cong \prod_{W_GC\backslash \Gen(C)} \bbQ\widehat{_p}[W_GC].$$
This induces a natural isomorphism
$$H^k_{W_GC}(C_GC\backslash X^C;\bbQ\widehat{_p}[\Gen(C)]) ~ \xrightarrow{\cong}
\prod_{W_GC\backslash \Gen(C)} H^k(C_GC\backslash X^C;\bbQ\widehat{_p}),$$
which comes from the adjunction $(i^*,i_!)$ of the functor restriction $i^*$ and coinduction $i_!$
for the ring homomorphism $i \colon \bbQ\widehat{_p} \to \bbQ\widehat{_p}[W_GC]$ and the obvious identification
$i_!(\bbQ\widehat{_p}) = \bbQ\widehat{_p}[W_GC]$.
There is an obvious bijection between the sets
$$\coprod_{(C) \in \calc_p} W_GC\backslash\Gen(C) \cong \con_p(G).$$
Now the claim follows from
Lemma~\ref{lem: rational computation of H^*_G(X;(bfK_{Bor})_{(0)})} and
Lemma~\ref{lem: identifying Q[Gen(C)] and theta_C cdot R(C) otimes_Z Q)}
\eqref{lem: identifying Q[Gen(C)] and theta_C cdot R(C) otimes_Z Q): Q}.
\end{proof}

\begin{theorem}[Computation of $K^n(EG \times_G X) \otimes_{\bbZ} \bbQ$]
\label{the: rational computation of K^*(EG times_GX)}
For every finite proper $G$-$CW$-complex $X$ and $n \in \bbZ$ there is a natural isomorphism
\begin{multline*}
\overline{\ch}_G^n \colon K^n(EG \times_G X) \otimes_{\bbZ} \bbQ
\\ \xrightarrow{\cong} ~
\prod_{i \in \bbZ} H^{2i+n}(G\backslash X;\bbQ) \times
\prod_{p} \prod_{(g) \in \con_p(G)}
H^{2i+n}(C_G\langle g \rangle\backslash X^{\langle g \rangle};\bbQ\widehat{_p}).
\end{multline*}
\end{theorem}
\begin{proof}
This follows from Lemma~\ref{lem: i^*_G(X;bfE) is bijective for finite X} and
Lemma~\ref{lem: refined non-multiplicative rational computation of H^*_G(X;(bfK_{Bor})_{(0)})}.
\end{proof}

\begin{lemma} \label{lem: comparison of G backslash Y and BG}
Let $Y \not= \emptyset$ be a proper $G$-$CW$-complex such that $\widetilde{H}_p(Y;\bbQ)$ vanishes for all $p$.
Let $f \colon Y \to \underline{E}G$ be a $G$-map.
Then $G\backslash f\colon G\backslash Y \to G\backslash \underline{E}G$ induces for all
$k$ isomorphisms
\begin{eqnarray*}
H_k(G\backslash f;\bbQ) \colon H_k(G\backslash Y;\bbQ) & \xrightarrow{\cong} &
H_k(G\backslash \underline{E}G;\bbQ);
\\
H^k(G\backslash f;\bbQ) \colon H^k(G\backslash \underline{E}G;\bbQ) & \xrightarrow{\cong} &
H^k(G\backslash Y;\bbQ);
\\
H^k(G\backslash f;\bbC) \colon H^k(G\backslash \underline{E}G;\bbC) & \xrightarrow{\cong} &
H^k(G\backslash Y,\bbC);
\\
H^k(G\backslash f;\bbQ\widehat{_p}) \colon H^k(G\backslash \underline{E}G;\bbQ\widehat{_p}) & \xrightarrow{\cong} &
H^k(G\backslash Y,\bbQ\widehat{_p});
\\
H^k(G\backslash f;\bbQ\widehat{_p} \otimes_{\bbQ} \bbC) \colon
H^k(G\backslash \underline{E}G;\bbQ\widehat{_p} \otimes_{\bbQ} \bbC) & \xrightarrow{\cong} &
H^k(G\backslash Y,\bbQ\widehat{_p} \otimes_{\bbQ} \bbC).
\end{eqnarray*}
\end{lemma}
\begin{proof}
The map $C_*(f) \otimes_{\bbZ} \id_{\bbQ} \colon C_*(G\backslash Y) \otimes_{\bbZ} \bbQ
\to C_*(\underline{E}G) \otimes_{\bbZ} \bbQ$ is  $\bbQ$-chain map of projective $\bbQ G$-chain
complexes and induces an isomorphism on homology. Hence it is  a $\bbQ G$-chain homotopy equivalence.
This implies that $C_*(f) \otimes_{\bbQ G} M$ and $\hom_{\bbQ G}(C_*(f),M)$ are
chain homotopy equivalences and induce isomorphisms on
homology and cohomology respectively for every $\bbQ$-module $M$.
\end{proof}

Now we can give the proof of Theorem~\ref{the: main theorem}.
\begin{proof}
We conclude from Lemma~\ref{lem: comparison of G backslash Y and BG}
that for any $g \in \con_p(G)$ the up to $C_G\langle g \rangle$-homotopy unique
$C_G\langle g \rangle$-map $f_g \colon EC_G\langle g \rangle \to \underline{E}C_G\langle g \rangle$
and the up to $G$-homotopy unique $G$-map $f \colon EG \to \underline{EG}$
induce isomorphisms
\begin{eqnarray}
H^k(G\backslash f;\bbQ) \colon H^k(G \backslash \underline{E}G;\bbQ)
& \xrightarrow{\cong} &
H^k(BG;\bbQ);
\label{identifying cohomology of underline{B}G and BG}
\\
\hspace{-5mm} H^k(C_G\langle g \rangle \backslash f_g;\bbQ\widehat{_p})
\colon H^k(C_G\langle g \rangle \backslash \underline{E}C_G\langle g \rangle;\bbQ\widehat{_p})
& \xrightarrow{\cong} &
H^k(BC_G\langle g \rangle,\bbQ\widehat{_p}).
\label{identifying cohomology of underline{B}C_G langle g rangle g and BC_Glangle g rangle}
\end{eqnarray}
Now apply Theorem~\ref{the: rational computation of K^*(EG times_GX)} to $X = \underline{E}G$
and use \eqref{identifying cohomology of underline{B}G and BG} and
\eqref{identifying cohomology of underline{B}C_G langle g rangle g and BC_Glangle g rangle}
together with the fact that $\underline{E}G^{\langle g \rangle}$ is a model for $\underline{E}C_G\langle g \rangle$.
\end{proof}

%%%%%%%%%%%%%%%%%%%%%%%%%%% Section 5 %%%%%%%%%%%%%%%%%%%%%%%%%%%%%%%%%%
\typeout{-----------------------  Section 5  ------------------------}

\tit{Multiplicative Structures}

In this section we want to deal with multiplicative structures and prove
Theorem~\ref{the: Multiplicative structure}.

\begin{remark}{\bf (Ring structures and multiplicative structures).}
\label{rem: ring spectra and Borel cohomology} \em
Suppose that the $\Omega$-spectrum $\bfE$ comes with the structure
of a ring spectrum $\mu \colon \bfE \wedge \bfE \to \bfE$. It induces
a multiplicative structure on the (non-equivariant) cohomology theory $H^*(-;\bfE)$ associated to $\bfE$.
Thus the equivariant cohomology theory
given by the equivariant Borel cohomology $H^*_?(E? \times_?-;\bfE)$
associated to $\bfE$ inherits a multiplicative structure the sense
of~\cite[Section~6]{Lueck(2004i)}.

If the contravariant $\GROUPOIDS$-$\Omega$-spectrum
$\bfF$ comes with a ring structure of  contravariant $\GROUPOIDS$-$\Omega$-spectra
$\mu \colon \bfF \wedge \bfF \to \bfF$, then the associated
equivariant cohomology theory $H^*_?(-;\bfF)$ inherits a
multiplicative structure. A ring structure on the $\Omega$-spectrum
$\bfE$ induces a ring structure of contravariant
$\GROUPOIDS$-$\Omega$-spectra on $\bfE_{\Bor}$. The induced
multiplicative structure on $H^*_?(-;\bfE_{\Bor})$ and the one on
$H^*_?(E? \times_?-;\bfE)$ are compatible with the natural
identification \eqref{H^n_G(X,A;bfE_{Bor})  = H^n(EG times_G (X,A);bfE)}.

A ring structure on the $\Omega$-spectrum $\bfE$ induces in a natural way a ring structure on
its rationalization $\Rat(\bfE)$. Thus a
ring structure on the contravariant $\GROUPOIDS$-$\Omega$-spectra on $\bfE_{\Bor}$
induces a
ring structure on the contravariant $\GROUPOIDS$-$\Omega$-spectra on $\left(\bfE_{\Bor}\right)_{(0)}$.
The natural transformation of equivariant
cohomology theories appearing in \eqref{i^*_?(-;bfE)} is compatible
with the induced  multiplicative structures.

In this discussion we are rather sloppy concerning the notion of a smash product.
Since we are not dealing with higher structures and just want to
take homotopy groups in the end, one can either use the classical approach
in the sense of Adams or the more advanced new constructions such as symmetric spectra.
\em
\end{remark}

\begin{lemma} \label{lem: multiplicative structures and Chern character}
The isomorphism appearing in
Lemma~\ref{lem: rational computation of H^*_G(X;(bfK_{Bor})_{(0)})}
is compatible with the multiplicative structure on the source and the one on the target
given by
$$(a,u_{p,(C)}) \cdot (b,v_{p,(C)}) ~ = ~ (a \cdot b,a \cdot v_{p,(C)´} + b \cdot v_{p,(C)} + u_{p,(C)} \cdot v_{p,(C)}),$$
for
\begin{eqnarray*}
(C) & \in & \calc_p(G);
\\
a,b & \in & H^*(BG;\bbQ);
\\
u_{p,(C)}, v_{p,(C)} & \in & H^{*}_ {W_GC}(C_GC \backslash X^C;\theta_C \cdot R(C) \otimes_{\bbZ} \bbQ\widehat{_p}),
\end{eqnarray*}
and the structures of a graded commutative ring on  $\prod_{i \in \bbZ} H^{2i+*}(BG;\bbQ)$ and
$\prod_{i \in \bbZ}H^{2i+*}_ {W_GC}(C_GC \backslash X^C;\theta_C \cdot R(C) \otimes_{\bbZ} \bbQ\widehat{_p})$
coming from the
cup-product and the multiplicative structure on $\theta_C \cdot R(C) \otimes_{\bbZ} \bbQ\widehat{_p}$
and the obvious $\prod_{i \in \bbZ} H^{2i+*}(BG;\bbQ)$-module structure
on $\prod_{i \in \bbZ}H^{2i+*}_ {W_GC}(C_GC \backslash X^C;\theta_C \cdot R(C) \otimes_{\bbZ} \bbQ\widehat{_p})$
coming from the canonical maps $C_GC\backslash X\to G\backslash X$ and $\bbQ \to \bbQ\widehat{_p}$.
\end{lemma}
\begin{proof} The proof consists of a straightforward calculation which is essentially based
on the following ingredients. In the sequel we use the notation of
\cite{Lueck(2004i)}.

The equivariant Chern character
of~\cite[Theorem~6.4]{Lueck(2004i)} is compatible with the multiplicative structures.

In Theorem~\ref{the: computation of K_*(BG) and K^*(BG) for finite G}
we have analyzed for every finite group $H$ the multiplicative structure on
$$
K^0(BH)  ~ \cong  ~ \bbZ \times \prod_{p}  \bbI_p(H)
\otimes_{\bbZ} \bbZ\widehat{_p}.
$$
Thus the Bredon cohomology group appearing in the target of the Chern character
whose source is $H_G^*\left(X;\left(\bfK_{\Bor}\right)_{(0)}\right)$ can be identified with
$$\left(\prod_{i \in \bbZ} H^{* + 2i}(G\backslash X;\bbQ)\right) \times \prod_p ~
\prod_{i \in \bbZ} ~ H^{2i +*}_{\bbQ\widehat{_p}\Sub(G;\calf)}(X;\bbI_p(?) \otimes_{\bbZ} \bbQ\widehat{_p})$$
with respect to the multiplicative structure analogously defined to the one appearing in
Theorem~\ref{the: computation of K_*(BG) and K^*(BG) for finite G} taking the obvious multiplicative structures
on the factors and the module structures of the factor for $p$ over
$\prod_{i \in \bbZ} H^{* + 2i}(G\backslash X;\bbQ)$ into account.

Fix a prime $p$. The $\bbQ\widehat{_p}[\aut(C)]$-map $R(C)\otimes_{\bbZ} \bbQ\widehat{_p} \to
\theta_C \cdot R(C)\otimes_{\bbZ} \bbQ\widehat{_p}$
given by multiplication with the idempotent
$\theta_C$ is compatible with the multiplicative structures.
Using the identification of Lemma~\ref{lem: T_C(I_p(C) otimes_Z Q_p) for C a finite cyclic p-group}
we obtain for each cyclic $p$-group $C$ a retraction compatible with the multiplicative structures.
$$\rho_C \colon \bbI_p(C) \otimes_{\bbZ} \bbQ\widehat{_p} ~ \to ~ T_C \left(\bbI_p(C) \otimes_{\bbZ} \bbQ\widehat{_p}\right)$$
Recall that $T_K\left(\bbI_p(K) \otimes_{\bbZ} \bbQ\widehat{_p}\right)$ is trivial unless $K$ is a non-trivial cyclic $p$-group.
Use these retractions as the maps $\rho_K$ in the definition of the isomorphism
$\nu$ of $\bbQ\widehat{_p}\Sub(G;\calf)$-modules for $M = \bbI_p(?) \otimes_{\bbZ} \bbQ\widehat{_p}$ in~\cite[(5.1)]{Lueck(2004i)}.
Then we obtain using the identification of Lemma~\ref{lem: T_C(I_p(C) otimes_Z Q_p) for C a finite cyclic p-group}
an isomorphism
of $\bbQ\widehat{_p}\Sub(G;\calf)$-modules
$$\bbI_p(?) \otimes_{\bbZ} \bbQ\widehat{_p} ~ \xrightarrow{\cong} ~
\prod_{(C) \in \calc_p} i(C)_!\left(\theta_C \cdot R(C) \otimes_{\bbZ} \bbQ\widehat{_p}\right),$$
which is compatible with the obvious multiplicative structure on the source and the one on the target given by
the product of the multiplicative structures on the factors
$i(C)_!\left(\theta_C \cdot R(C) \otimes_{\bbZ} \bbQ\widehat{_p}\right)$ coming from the obvious one on
$\theta_C \cdot R(C) \otimes_{\bbZ} \bbQ\widehat{_p}$. Using the adjunction $(i(C)^*,i(C)_!)$ this isomorphism
induces an $\bbQ\widehat{_p}$-isomorphism compatible with the multiplicative structures
$$H^n_{\bbQ\widehat{_p}\Sub(G;\calf)}(X;\bbI_p(?) \otimes_{\bbZ} \bbQ\widehat{_p}) ~ \xrightarrow{\cong} ~
\prod_{(C) \in \calc_p(G)} H^n_{W_GC}(C_GC\backslash X;\theta_C \cdot R(C) \otimes_{\bbZ} \bbQ\widehat{_p}).$$

\end{proof}

Because the isomorphism in Lemma~\ref{lem: identifying Q[Gen(C)] and theta_C cdot R(C) otimes_Z Q)}
\eqref{lem: identifying Q[Gen(C)] and theta_C cdot R(C) otimes_Z Q): C}
is compatible with the multiplicative structures, it implies together with
Lemma~\ref{lem: multiplicative structures and Chern character}

\begin{lemma} \label{lem: complex computation of H^*_G(X;(bfK_{Bor})_{(0)})}
For every proper $G$-$CW$-complex $X$ and $n \in \bbZ$ there is a $\bbC$-isomorphism,
natural in $X$,
\begin{multline*}
\overline{\ch}^n_{G,\bbC} \colon H^*_G(X;(\bfK_{\Bor})_{(0)})\otimes_{\bbQ} \bbC\xrightarrow{\cong}
\\
\left(\prod_{i \in \bbZ} H^{2i+n}(G\backslash X;\bbC)\right) \times
\prod_{p} ~ \prod_{(g) \in \con_p(G)}
\left(\prod_{i \in \bbZ} H^{2i+n}(C_G\langle g \rangle\backslash X^{\langle g \rangle};\bbQ\widehat{_p} \otimes_{\bbQ} \bbC)\right),
\end{multline*}
which is compatible with the  multiplicative structure on the target given by
$$
\left(a, u_{p,(g)}\right) \cdot \left(b, v_{p,(g)}\right)
~ = ~
\left(a\cdot b, (a \cdot v_{p,(g)} + b \cdot u_{p,(g)} + u_{p,(g)} \cdot v_{p,(g)})\right)
$$
for
\begin{eqnarray*}
(g) & \in & \con_p(G);
\\
a,b  & \in & \prod_{i \in \bbZ} H^{2i+*}(G\backslash X;\bbC),
\\
u_{p,(g)}, v_{p,(g)} & \in & \prod_{i \in \bbZ}
H^{2i+*}(C_G\langle g \rangle\backslash X^{\langle g \rangle};\bbQ\widehat{_p} \otimes_{\bbQ} \bbC),
\end{eqnarray*}
and the structures of a graded commutative ring on  $\prod_{i \in \bbZ} H^{2i+*}(G\backslash X;\bbC)$ and
$\prod_{i \in \bbZ} H^{2i+*}(C_G\langle g \rangle\backslash X^{\langle g \rangle};\bbQ\widehat{_p} \otimes_{\bbQ} \bbC)$
coming from the cup-product and the obvious $\prod_{i \in \bbZ} H^{2i+*}(G\backslash X;\bbC)$-module structure on
$\prod_{i \in \bbZ} H^{2i+*}(C_G\langle g \rangle\backslash X^{\langle g \rangle};\bbQ\widehat{_p} \otimes_{\bbQ} \bbC)$
coming from the canonical map $BC_G\langle g \rangle \to BG$.
\end{lemma}

Now we are ready to prove Theorem~\ref{the: Multiplicative structure}
\begin{proof}
The isomorphism appearing in Lemma~\ref{lem: i^*_G(X;bfE) is bijective for finite X}
is compatible with the multiplicative structures. This is also true for the versions of isomorphisms
\eqref{identifying cohomology of underline{B}G and BG} and
\eqref{identifying cohomology of underline{B}C_G langle g rangle g and BC_Glangle g rangle},
where the coefficients $\bbQ$ and $\bbQ\widehat{_p}$ are replaced by $\bbC$ and $\bbQ\widehat{_p} \otimes_{\bbQ} \bbC$.
Now put these together with the isomorphism appearing in
Lemma~\ref{lem: complex computation of H^*_G(X;(bfK_{Bor})_{(0)})}.
\end{proof}

\begin{remark}{\bf (Difference between rationalization and com\-plexi\-fi\-ca\-tion).}
\label{rem: Difference between rationalization and
complexification} \em First of all we want to emphasize that the
isomorphism appearing in Theorem~\ref{the: Multiplicative
structure} is \emph{not} obtained from the isomorphism appearing
in Theorem~\ref{the: main theorem} by applying $ - \otimes_{\bbQ}
\bbC$ since the corresponding statement is already false for the
two isomorphisms appearing in Lemma~\ref{lem: identifying
Q[Gen(C)] and theta_C cdot R(C) otimes_Z Q)}. Moreover, the
isomorphism appearing in Theorem~\ref{the: main theorem} is
\emph{not} compatible with the standard multiplicative structures
on the source and the multiplicative structure on the target which
is defined analogously  to the one on the complexified target in
Theorem~\ref{the: Multiplicative structure}. The reason is that
the isomorphism appearing in Lemma~\ref{lem: identifying Q[Gen(C)]
and theta_C cdot R(C) otimes_Z Q)}~ \eqref{lem: identifying
Q[Gen(C)] and theta_C cdot R(C) otimes_Z Q): Q} cannot be chosen
to be compatible with the obvious multiplicative structure on its
target if we use on the source the multiplicative structure coming
from the obvious identification $\bbQ\widehat{_p}[\Gen(C)] =
\prod_{\Gen(C)} \bbQ\widehat{_p}$ and the product ring structure
on $\prod_{\Gen(C)} \bbQ\widehat{_p}$.

One can easily check by hand that there is \emph{no}
 $\bbQ\widehat{_3}[\aut(\bbZ/3)]$-isomorphism compatible with the multiplicative structures
$$\theta_{\bbZ/3} \cdot R(\bbZ/3)  \otimes \bbQ\widehat{_3}
= \bbI_{\bbZ/3} \otimes \bbQ\widehat{_3}  ~ \xrightarrow{\cong} ~
\bbQ\widehat{_3} \times \bbQ\widehat{_3}$$ if we equip the target
with the $\aut(\bbZ/3) \cong \bbZ/2$-action given by flipping the
factors and the product $\bbQ$-algebra structure. The point is
that $\bbQ\widehat{_3}$ does not contain a primite $3$-rd root of
unity (in contrast to $\bbC$, see Lemma~\ref{lem: identifying
Q[Gen(C)] and theta_C cdot R(C) otimes_Z Q)}~\eqref{lem:
identifying Q[Gen(C)] and theta_C cdot R(C) otimes_Z Q): C} ). \em
\end{remark}

\begin{example}[Multiplicative structure over $\bbQ$]
  \label{exa: including multiplicative structures} \em
In general we can give a simple formula for the multiplicative structure only after complexifying
as explained in Remark~\ref{rem: Difference between rationalization and complexification}.
In the following special case this can be done already after rationalization.
Suppose that for any non-trivial cyclic subgroup $C$ of prime power order
$\widetilde{H}^n(BC_GC;\bbQ) = 0$ holds for all $n \in \bbZ$ and that $W_GC = \aut(C)$.
The latter means that any automorphism of $C$
is given by conjugation with some element in $N_GC$. Suppose furthermore that there is a
finite model for $\underline{E}G$. Then we obtain $\bbQ$-isomorphisms
\begin{eqnarray*}
K^0(BG) \otimes_{\bbZ} \bbQ &  \cong &
\prod_{i \in \bbZ} H^{2i}(BG;\bbQ) \times \prod_p (\bbQ\widehat{_p})^{r_p(G)};
\\
K^1(BG) \otimes_{\bbZ} \bbQ &  \cong & \prod_{i \in \bbZ} H^{2i+1}(BG;\bbQ),
\end{eqnarray*}
where $r_p(G)$ is the number of conjugacy classes of non-trivial cyclic subgroups
of $p$-power order what is in this situation the same as $|\con_p(G)|$.
The isomorphisms above are compatible with the multiplicative structure on the target given by
\begin{eqnarray*}
(a,u) \cdot (b,v) & = & (a \cup  b, a_0 \cdot v + b_0 \cdot u + u \cdot v);
\\
(a,u) \cdot c & = & a \cup c;
\\
c \cdot d & = & c \cup d,
\end{eqnarray*}
for $a,b \in \prod_{i \in \bbZ} H^{2i}(BG;\bbQ)$, $c,d \in \prod_{i \in \bbZ} H^{2i+1}(BG;\bbQ)$
and $u,v \in  \prod_p (\bbQ\widehat{_p})^{r_p(G)}$, where $a_0 \in \bbQ$ and $b_0 \in \bbQ$ are the components
of $a$ and $b$ in $H^0(BG;\bbQ) = \bbQ \cdot 1$ and we equip
$\prod_p (\bbQ\widehat{_p})^{r_p(G))}$ with the structure of a $\bbQ$-algebra
coming from the product of the obvious $\bbQ$-algebra structures on the various factors $\bbQ\widehat{_p}$.
This follows from Lemma~\ref{lem: multiplicative structures and Chern character} and
the conclusion from the formula $\theta_C \cdot \theta_C = \theta_C$ and
Lemma~\ref{lem: identifying Q[Gen(C)] and theta_C cdot R(C) otimes_Z Q)}
\eqref{lem: identifying Q[Gen(C)] and theta_C cdot R(C) otimes_Z Q): Q}
that $\left(\theta_C \cdot R(C)\right)^{\aut(C)}$
is generated as $\bbQ$-vector space by $\theta_C$ and hence is as $\bbQ$-algebra
isomorphic to $\bbQ$.

If we  furthermore assume that $\widetilde{H}_n(BG;\bbQ) = 0$ for all $n \in \bbZ$,
the formula simplifies to
\begin{eqnarray*}
K^0(BG) \otimes_{\bbZ} \bbQ &  \cong &
\bbQ \times \prod_p (\bbQ\widehat{_p})^{r_p(G)};
\\
K^1(BG) \otimes_{\bbZ} \bbQ &  \cong & 0.
\end{eqnarray*}
The first isomorphism is compatible with the multiplicative
structures if we put on the target the one given by
$$(m,a) \cdot (m,b) ~ = ~ (mn,m\cdot b + n \cdot a + a \cdot b)$$
for $m,n \in \bbQ$, $a,b \in \prod_p (\bbQ\widehat{_p})^{r_p(G)}$
and we equip $\prod_p (\bbQ\widehat{_p})^{r_p(G)}$ with the structure of a $\bbQ$-algebra
coming from the product of the obvious $\bbQ$-algebra structures on the various factors $\bbQ\widehat{_p}$.

\em
\end{example}

%%%%%%%%%%%%%%%%%%%%%%%%%%% Section 6 %%%%%%%%%%%%%%%%%%%%%%%%%%%%%%%%%%
\typeout{-----------------------  Section 6  ------------------------}

\tit{Weakening the Finiteness Conditions}

In this section we want to weaken the finiteness assumption occurring
in Theorem~\ref{the: main theorem} and
Theorem~\ref{the: Multiplicative structure}.

A $\bbZ$-module $M$ is \emph{almost trivial}
if there is an element $r \in \bbZ, r \not=0$ such that $rm = 0$ holds for
all $m \in M$. A $\bbZ$-module $M$ is \emph{almost finitely generated} if
$M/\tors(M)$ is a finitely generated $\bbZ$-module and $\tors(M)$ is almost trivial.
A $\bbZ$-homomorphism is an \emph{almost isomorphism} if its kernel and cokernel are
almost trivial. An almost isomorphism becomes an isomorphism after
rationalization.

The full subcategories of the category of $\bbZ$-modules given by almost trivial submodules and by almost finitely
generated submodules are Serre-subcategories, i.e. are closed under subobjects, quotients, and extensions.
In particular there is a Five-Lemma for almost isomorphisms. These notions and facts are introduced and
proved in~\cite[Section~4]{Lueck-Reich-Varisco(2003)}.

The main result of this section is:

\begin{theorem}[Weakening the finiteness assumption]
\label{the: Weakening the finiteness assumption} The conclusions of
Theorem~\ref{the: main theorem} and
Theorem~\ref{the: Multiplicative structure} remain true if we replace the condition
that there is a
cocompact  $G$-$CW$-model for the classifying space $\underline{E}G$ for
proper $G$-actions by the following weaker set of conditions:

There exists a $G$-$CW$-complex $X$ satisfying:
\begin{enumerate}

\item \label{the: Weakening the finiteness assumption: condition 1}
The $G$-$CW$-complex $X$ is proper and finite dimensional.
There is an upper bound on the orders of its isotropy groups.
The set of conjugacy classes $(C)$ of finite cyclic subgroups $C \subseteq G$
of prime power order with $X^C \not= \emptyset$ is finite;

\item \label{the: Weakening the finiteness assumption: condition 2}
For all for $k \in \bbZ$ we have $H_k(X;\bbZ) \cong H_k(\pt;\bbZ)$;

\item \label{the: Weakening the finiteness assumption: condition 3}
For any finite cyclic subgroup of prime power order $C \subseteq G$ and integer $k$ the $\bbZ$-module
$H_k(X^C;\bbZ)$ is almost finitely generated;

\item \label{the: Weakening the finiteness assumption: condition 4}
For any finite cyclic subgroup of prime power order $C \subseteq G$ and integer $k$ the $\bbZ$-module
$H_k(C_GC\backslash X^C;\bbZ)$ is almost finitely generated.

\end{enumerate}

If $X$ satisfies conditions \eqref{the: Weakening the finiteness assumption: condition 1},
\eqref{the: Weakening the finiteness assumption: condition 2} and
\eqref{the: Weakening the finiteness assumption: condition 3} above, then the condition
\eqref{the: Weakening the finiteness assumption: condition 4} is satisfied if and only
if for any finite cyclic subgroup of prime power order $C \subseteq G$ and integer $k$ the $\bbZ$-module
$H_k(BC_GC;\bbZ)$ is almost finitely generated.
\end{theorem}

\begin{remark}[Weakening the finiteness conditions for $\underline{E}G$]
\label{rem: discussing the finiteness assumptions} \em
Notice that the conditions  \eqref{the: Weakening the finiteness assumption: condition 1},
\eqref{the: Weakening the finiteness assumption: condition 2},
\eqref{the: Weakening the finiteness assumption: condition 3} and
\eqref {the: Weakening the finiteness assumption: condition 4} in
Theorem~\ref{the: Weakening the finiteness assumption} are satisfied, if
the set of conjugacy classes of finite subgroups of $G$ is finite,
there is a finite dimensional model for $\underline{E}G$ and
for any finite cyclic subgroup of prime power order $C \subseteq G$ and integer $k$ the $\bbZ$-module
$H_k(BC_GC;\bbZ)$ is almost finitely generated.
\em
\end{remark}

\begin{remark}[Virtually torsionfree groups] \label{rem: type of underline EG} \em
Suppose that $G$ contains a torsionfree subgroup $H \subseteq G$ of finite index.
If there is a finite dimensional model for $BH$, then there exists a
finite dimensional model for $\underline{E}G$ ~\cite{Serre(1971)}).
However, if there is a finite model for $BH$,
this does not implies that $G$ has only finitely many conjugacy classes of subgroups
or that there is a cocompact model for $\underline{E}G$ or that the centralizers
$C_GC$ of finite cyclic subgroups are finitely generated
\cite[Section~7]{Leary-Nucinkis(2003)}.
\em
\end{remark}

The proof of Theorem~\ref{the: Weakening the finiteness assumption} needs some preparation.

\begin{lemma} \label{lem: almost iso induced by EG times_G X to X/G}
Let $X$ be a proper $G$-$CW$-complex. Let $\pr\colon EG \times_G X \to G\backslash X$ be the projection.
Fix an integer $n \in \bbZ$.
\begin{enumerate}
\item \label{lem: almost iso induced by EG times_G X to X/G: almost isomorphism}
Suppose that there exists for
each $m \ge 0$ a   positive integer $d(m)$ such that for any isotropy group
$H$ of $X$ multiplication with $d(m)$ annihilates $\widetilde{H}_m(BH;\bbZ)$.
Then the induced map
$$H_n(\pr;\bbZ) \colon H_n(EG \times_G X;\bbZ) \to  H_n(G\backslash X;\bbZ).$$
is an almost isomorphism for all $n \in \bbZ$;

\item \label{lem: almost iso induced by EG times_G X to X/G: rational isomorphism}
The induced map
$$H_n(\pr;\bbQ) \colon H_n(EG \times_G X;\bbQ) \to  H_n(G\backslash X;\bbQ).$$
is a $\bbQ$-isomorphism.
\end{enumerate}
\end{lemma}
\begin{proof}
\eqref{lem: almost iso induced by EG times_G X to X/G: almost isomorphism}
This is proved in~\cite[Lemma~8.1]{Lueck-Reich-Varisco(2003)}.
\\
\eqref{lem: almost iso induced by EG times_G X to X/G: rational isomorphism}
The proof is analogous to the one of \eqref{lem: almost iso induced by EG times_G X to X/G: rational isomorphism}.
\end{proof}

The next result is a generalization of Lemma~\ref{lem: i^*_G(X;bfE) is bijective for finite X} in
the case $\bfE = \bfK$.

\begin{lemma}
\label{lem: i^*_G(X;bfK) is bijective for certain X}
Let $X$ be  a finite dimensional proper $G$-$CW$-complex such that there is a
bound on the orders of finite subgroups.
Let $J$ be  the set of conjugacy classes $(C)$ of finite cyclic subgroups $C \subseteq G$
of prime power order with $X^C \not= \emptyset$. Suppose that $|J|$ is finite.
Furthermore assume that $H_k(C_GC\backslash X^C;\bbZ)$ is almost finitely generated for
every $k \in \bbZ$ and every finite cyclic subgroup of prime power order $C \subseteq G$.

Then the map
$$i^n_G(X;\bfK) \colon H_G^n\left(X;\bfK_{\Bor}\right) \otimes_{\bbZ} \bbQ \to
H_G^n\left(X;\left(\bfK_{\Bor}\right)_{(0)}\right)$$
is bijective for all $n  \in \bbZ$.
\end{lemma}
\begin{proof} In the sequel we use the notation of ~\cite{Lueck(2004i)}.
Let $\calf(X)$ be the set of conjugacy classes
of subgroups of $H \subseteq G$ with $X^H \not= \emptyset$.
Since $X$ is proper and has finite orbit type, $\calf(X)$ is finite and
$(H) \in \calf(X)$ implies that $H$ is finite.
Since $X$ is proper and finite dimensional,
there is a spectral sequence converging to  $H_G^{s+t}\left(X;\bfK_{\Bor}\right)$
whose $E_2$-term is $E_2^{s,t} = H^s_{\Sub(G;\calf(X))}(X;K^t(BH))$ and a
spectral sequence converging to  $H_G^{s+t}\left(X;\left(\bfK_{\Bor}\right)_{(0)}\right)$
whose $E_2$-term is $E_2^{s,t} = H^s_{\Sub(G;\calf(X))}(X;K^t(BH)\otimes_{\bbZ} \bbQ)$.
Since $\bbQ$ is flat over $\bbZ$, it suffices to show that the canonical map
$$H^s_{\Sub(G;\calf(X))}(X;K^t(BH))\otimes_{\bbZ} \bbQ  \to
H^s_{\Sub(G;\calf(X))}(X;K^t(BH)\otimes_{\bbZ} \bbQ)$$
is bijective for all $s$ and $t$.
We have already explained that the contravariant $\bbZ\Sub(G;\calf(X))$-module sending
$H$ to $K^t(BH)$ is zero for odd $t$  and given for even $t$ by
$$K^t(BH) ~ \cong ~ \bbZ \times \prod_{p} \bbI_p(H) \otimes_{\bbZ} \bbZ\widehat{_p}.$$
One easily checks using Lemma~\ref{lem: T_H(I_p(H) otimes_Z Q_p) for H non finite cyclic p-group}
and \eqref{computation of T_C (bbI_p(C)) for C Z/p^k} for any finite group $H$
$$T_H\left(K^0(BH)\right)
 \cong  \left\{
\begin{array}{lll}
\bbZ & & H = \{1\};
\\
\ker\left(\res_H^{H'} \colon R(H) \to R(H')\right) \otimes_{\bbZ} \bbZ\widehat{_p} & &
H \text{ cyclic } p\text{-group},
\\
& & H' \subseteq H, [H:H'] = p;
\\
0& & \text{ otherwise}.
\end{array}\right.
$$
For every finite cyclic subgroup $K \subseteq G$ of order $p^r$ for some prime $p$ and integer $r \ge 1$
choose a retraction $r'(K) \colon K^0(BK) \to T_K( K^0(BK))$ of the $\bbZ$-homomorphism
$j(K)\colon T_K( K^0(BK)) \to K^0(BK)$ given by inclusion. Such $r'(K)$ exists
since  $R(K')$ and hence the image of $\res_K^{K'} \colon R(K) \to R(K')$ is a finitely generated free $\bbZ$-module
what implies that
$\ker\left(\res_K^{K'} \colon R(K) \to R(K')\right)$ is a direct summand
of the finitely generated free $\bbZ$-module $\bbI_p(K) = \bbI(K)$.
Since $W_GK$ is finite, we can define a
$\bbZ[W_GK]$-map
$$r(K) \colon K^0(BK) \to T_K( K^0(BK)),
\quad x \mapsto \sum_{g \in W_GK} g \cdot r'(K)(g^{-1} \cdot x).$$
Then $r(K) \circ j(K) = |W_GK| \cdot \id$.
For $K = \{1\}$ let $r(K) \colon K^0(BK) \xrightarrow{\cong} \bbZ$ be the obvious isomorphism
which is for trivial reasons a $W_GK$-map.
Define a map of contravariant
$\Sub(G;\calf(X))$-modules
$$\nu \colon K^0(B?) \to \prod_{(K) \in J} i(K)_!T_K(K^0(B?))$$
by requiring that the composite of $\nu$ with the projection onto the factor belonging to $(K) \in J$ is the adjoint
for the pair $(i(K)^*,i(K)_!)$ of the $W_GK$-map $r(K)$. Analogously to
the proof of~\cite[Theorem~2.14~(b)]{Lueck(2004i)}
one shows that $\nu(H)$
is injective for all objects $H \in \Sub(G;\calf(X))$. Here we use the fact that
$r(K) \circ j(K)$ is injective for all $K$ with $(K) \in J$
since $r(K) \circ j(K) = |W_GK| \cdot \id$ and $ K^0(BK)$ and hence $T_K( K^0(BK))$ is torsionfree. Then one constructs
analogously to  the proof of~\cite[Theorem~5.2]{Lueck(2004i)}
for each object $(H) \in \Sub(G;\calf(X))$ a $\bbZ$-homomorphism
$$\mu(H) \colon \left(\prod_{(K) \in J} i(K)_!T_K(K^0(B?))\right)(H) ~ \to ~ K^0(BH)$$
and checks that $\nu(H) \circ \mu(H)$ can be written as a diagonal matrix $A(H)$ which has upper
triangular form and has maps of the shape $r \cdot\id$ as diagional entry,
where each $r$ divides a certain integer  $M(|H|)$ depending only on the order
of $|H|$.  There is an integer $N(|H|)$ depending
only on the order of $|H|$ such that the size of the square matrix $A(H)$ is bounded by $N(|H|)$.
The existence of the numbers $M(|H|)$ and  $N(|H|)$ follows from the finiteness of
$J$.  Hence for each object $H \in \Sub(G;\calf(X))$ the cokernel of $\nu(H)$ is annihilated by $M(|H|)^{N(|H|)}$.
Since there is an upper bound on the orders of finite subgroups of $G$, we can find an integer $L$
such that for each object $H \in \Sub(G;\calf(X))$ the cokernel of $\nu(H)$ is annihilated by $L$.
The short exact sequence of $\bbZ\Sub(G;\calf(X))$-modules
 $$0 \to K^0(B?) \xrightarrow{\nu} \prod_{(K) \in J} i(K)_!T_K(K^0(B?))
\xrightarrow{\pr} \cok(\nu) \to 0$$
induces a long exact sequence
\begin{multline*}
\ldots \to H^{s-1}_{\Sub(G;\calf(X))}(X;\cok(\nu)) \to
 H^{s}_{\Sub(G;\calf(X))}(X;K^0(B?))
\\\to
H^{s}_{\Sub(G;\calf(X))}(X;\prod_{(K) \in J} i(K)_!T_K(K^0(B?)))
\to H^{s}_{\Sub(G;\calf(X))}(X;\cok(\nu)) \to \ldots
\end{multline*}
Since multiplication with $L$ induces the zero map $\cok(\nu) \to \cok(\nu)$,
multiplication with $L$ induces also the zero map on
$H^{s}_{\Sub(G;\calf(X))}(X;\cok(\nu))$. Hence
$H^{s}_{\Sub(G;\calf(X))}(X;\cok(\nu)) \otimes_{\bbZ} \bbQ$ is trivial. Since $J$ is finite,
and $- \otimes_{\bbZ} \bbQ$ is an exact functor which commutes with finite products,
we obtain from the adjunction $(i(K)^*.i(K)_!)$ a natural isomorphism
$$H^{s}_{\Sub(G;\calf(X))}(X;K^0(B?)) \otimes_{\bbZ} \bbQ
~ \xrightarrow{\cong} ~
\prod_{(K) \in J} H^s_{W_GK}(C_GK\backslash X^K;T_KK^0(B?))  \otimes_{\bbZ} \bbQ .$$
Similarly we get an isomorphism
$$H^{s}_{\Sub(G;\calf(X))}(X;K^0(B?)\otimes_{\bbZ} \bbQ)
~ \xrightarrow{\cong} ~
\prod_{(K) \in J} H^s_{W_GK}(C_GK\backslash X^K;T_KK^0(B?)  \otimes_{\bbZ} \bbQ) .$$
Hence it remains to show for each $(K) \in J$ and $s \ge 0$ that the canonical map
$$
H^s_{W_GK}(C_GK\backslash X^K;T_KK^0(B?))  \otimes_{\bbZ} \bbQ
~ \to ~
H^s_{W_GK}(C_GK\backslash X^K;T_KK^0(B?)  \otimes_{\bbZ} \bbQ)
$$
is bijective. Abbreviate $C_* = C_*(C_GK\backslash Y^K)$, $L = W_GK$ and $M = T_KK^0(B?)$.
Then $L$ is a finite group, $C_*$ is a  $\bbZ L$-chain complex
which is free over $\bbZ$ and for which there exists an integer
$n \ge 1$ such that $\tors(H_s(C_*))$ is annihilated by $n$ and
$H_s(C_*)/\tors(H_s(C_*))$ is a finitely generated $\bbZ$-module for all $s$. It remains to show for
the $\bbZ L$-module $M$ that the canonical map
$$H^s(\hom_{\bbZ L}(C_*,M))  \otimes_{\bbZ} \bbQ ~ \to ~
H^s(\hom_{\bbZ L}(C_*,M \otimes_{\bbZ} \bbQ))$$ is bijective for
all $s$. Let $i \colon \{1\} \to L$ be the inclusion. Since $L$ is
finite, induction $i_*$ and coinduction $i_!$ agree. Hence we get
natural $\bbZ L$-chain maps $a_* \colon C_* \to  i_*i^*C_*$ and
$b_*\colon i_*i^*C_* \to C_*$ such that $b_* \circ a_*$ is
multiplication with $|L|$. They are explicitly given by
\begin{eqnarray*}
a_s \colon  C_s \to \bbZ L \otimes_{\bbZ} C_s, & \quad &
x \mapsto \sum_{l \in L} l \otimes l^{-1} \cdot x;
\\
b_s \colon  \bbZ L \otimes_{\bbZ} C_s  \to  C_s, & \quad &
l \otimes y \mapsto l \cdot y.
\end{eqnarray*}
Hence we obtain a commutative diagram
$$
\begin{CD}
H^s(\hom_{\bbZ L}(C_*,M)) \otimes_{\bbZ} \bbQ
@>>>
H^s(\hom_{\bbZ L}(C_*,M \otimes_{\bbZ} \bbQ))
\\
@V H^s(b_*) \otimes_{\bbZ} \id_{\bbQ} VV @V H^s(b_*)  VV
\\
H^s(\hom_{\bbZ L}(\bbZ L \otimes_{\bbZ} C_*,M)) \otimes_{\bbZ} \bbQ
@>>>
H^s(\hom_{\bbZ L}(\bbZ L \otimes_{\bbZ} C_*,M\otimes_{\bbZ} \bbQ))
\\
@V H^s(a_*) \otimes_{\bbZ} \id_{\bbQ} VV @V H^s(a_*) VV
\\
H^s(\hom_{\bbZ L}(C_*,M)) \otimes_{\bbZ} \bbQ
@>>>
H^s(\hom_{\bbZ L}(C_*,M \otimes_{\bbZ} \bbQ))
\end{CD}
$$
where the horizontal arrows are the canonical maps  and the composite of the two left vertical maps
and the composite of the two right vertical maps are isomorphisms. Hence it suffices to show
that the middle horizontal  arrow is an isomorphism. It can be identified with the canonical map
$$
H^s(\hom_{\bbZ}( C_*,M)) \otimes_{\bbZ} \bbQ
~ \to ~
H^s(\hom_{\bbZ}(C_*,M \otimes_{\bbZ} \bbQ)).$$
Notice for the sequel that $C_*$ is $\bbZ$-free. By the universal coefficient theorem
we get a commutative diagram with exact rows and the canonical maps as horizontal arrows
$$
\begin{CD}
0 & &  0
\\
@VVV @VVV
\\
\ext_{\bbZ}(H_{s-1}(C_*),M) \otimes_{\bbZ} \bbQ
@>>>
\ext_{\bbZ}(H_{s-1}(C_*),M \otimes_{\bbZ} \bbQ)
\\
@VVV @VVV
\\
H^s(\hom_{\bbZ}(C_*,M)) \otimes_{\bbZ} \bbQ
@>>>
H^s(\hom_{\bbZ}(C_*,M \otimes_{\bbZ} \bbQ))
\\
@VVV @VVV
\\
\hom_{\bbZ}(H_s(Z),M) \otimes_{\bbZ}\bbQ
@>>>
\hom_{\bbZ}(H_s(Z),M \otimes_{\bbZ} \bbQ))
\\
@VVV @VVV
\\
0 & & 0
\end{CD}
$$

Since $H_s(C_*)/\tors(H_s(C_*)$
is a finitely generated free abelian group and there is an integer $n$ which annihilates
$\tors(H_s(C_*))$, the rational vector spaces
$\ext_{\bbZ}(H_{s-1}(C_*),M) \otimes_{\bbZ} \bbQ$ and $\ext_{\bbZ}(H_{s-1}(C_*),M \otimes_{\bbZ} \bbQ)$
vanish and the lower vertical arrow is bijective. Hence the middle arrow is bijective.
This finishes the proof of Lemma~\ref{lem: i^*_G(X;bfK) is bijective
  for certain X}.
\end{proof}

\begin{lemma} \label{H^*,-,(K_Bor)_(0))-iso}
Let $X$ be a proper $G$-$CW$-complex such that for any
cyclic subgroup $C \subseteq G$ of prime power order and any $k \in \bbZ$ we have $H_k(X^C;\bbQ) \cong H_k(\pt;\bbQ) $.
Then the up to $G$-homotopy
unique $G$-map $f \colon X \to \underline{E}G$ induces for every $n \in \bbZ$ an isomorphism
$$H_G^n\left(f;\left(\bfK_{\Bor}\right)_{(0)}\right) \colon
H_G^n\left(\underline{E}G;\left(\bfK_{\Bor}\right)_{(0)}\right)
\xrightarrow{\cong}
H_G^n\left(X;\left(\bfK_{\Bor}\right)_{(0)}\right).$$
\end{lemma}
\begin{proof}
Because of Lemma~\ref{lem: refined non-multiplicative rational computation of H^*_G(X;(bfK_{Bor})_{(0)})}
it suffices to show for every $n \in \bbZ$ and every cyclic subgroup of prime power order that the map
$$H^n(C_GC\backslash f^C;M) \colon H^n(C_GC\backslash (\underline{E}G)^C;M) \to H^n(C_GC\backslash X^C;M)$$
is bijective for any $\bbQ$-module $M$. The Atiyah-Hirzebruch spectral sequence for the fibration
$X^C \to EC_GC \times_{C_GC} X^C \to BC_GC$ together with the vanishing of
$\widetilde{H}_*(X^C;\bbQ)$ implies that the projection
$\pr \colon EC_GC \times_{C_GC} X^C \to BC_GC$ induces for all $n \in \bbZ$ isomorphisms
$$H_n(\pr;\bbQ) \colon H_n(EC_GC \times_{C_GC} X^C;\bbQ) \to H_n(BC_GC;\bbQ).$$
The projection $\pr' \colon EC_GC \times_{C_GC} X^C \to C_GC \backslash X^C$ induces for all
$n \in \bbZ$ isomorphisms
$$H_n(\pr';\bbQ) \colon H_n(EC_GC \times_{C_GC} X^C ;\bbQ) \to H_n(C_GC\backslash X^C;\bbQ)$$
by Lemma~\ref{lem: almost iso induced by EG times_G X to X/G}~
\eqref{lem: almost iso induced by EG times_G X to X/G: rational isomorphism}.
The same is also true for $\underline{E}G$ instead of $X$.
This implies that
$$H_n(C_GC\backslash f^C;\bbQ) \colon H_n(C_GC\backslash X^C;\bbQ) \to H_n(C_GC\backslash  (\underline{E}G)^C;\bbQ)$$
is bijective for all $n \in \bbZ$. Hence
$H^n(C_GC\backslash f^C;M)$ is bijective for any $\bbQ$-module $M$.
\end{proof}

\begin{lemma} \label{lem: Smith theory}
Let $X$ be a proper finite dimensional  $G$-$CW$-complex such that $H_k(X;\bbZ) \cong H_k(\pt;\bbZ)$ holds for any $k \in \bbZ$.
Let $C \subseteq G$ be a cyclic group of prime power order. Suppose that
$\widetilde{H}_n(X^C;\bbZ)$ is almost finitely generated for each $n \in \bbZ$.
Then:

\begin{enumerate}

\item \label{lem: Smith theory: almost trivial}
The $\bbZ$-module $\widetilde{H}_n(X^C;\bbZ)$ is almost trivial and the $\bbQ$-module $\widetilde{H}_n(X^C;\bbQ)$
is trivial for all $n \in \bbZ$;

\item \label{lem: Smith theory: almost isomorphism}
The map $H_n(\pr;\bbZ) \colon H_n(EC_GC \times_{C_GC} X^C;\bbZ) \to H_n(BC_GC;\bbZ)$ induced by the
projection $\pr \colon EC_GC \times_{C_GC} X^C \to BC_GC$ is an almost isomorphism for
all $n \in \bbZ$.

%\item \label{lem: Smith theory: rational isomorphism}
%The map $H_n(\pr;\bbQ) \to H_n(EC_GC \times_{C_GC} X^C;\bbQ) \to H_n(BC_GC;\bbQ)$ is an isomorphism for
%all $n \in \bbZ$.

\end{enumerate}
\end{lemma}
\begin{proof}
\eqref{lem: Smith theory: almost trivial}
Suppose that $C$ has order $p^k$ for $k \ge 1$. Then
$\widetilde{H}_n(X;\bbF_p) = 0$ for all $n \in \bbZ$ if $\bbF_p$ is the finite field of
order $p$. By Smith theory $\widetilde{H}_n(X^C;\bbF_p) = 0$ for all $n \in \bbZ$
\cite[Theorem~5.2]{Bredon(1972)}. This implies by the Bockstein sequence
associated to $0 \to \bbZ \xrightarrow{p \cdot \id} \bbZ \to \bbF_p \to 0$
that $p \cdot \id \colon \widetilde{H}_n(X^C;\bbZ) \to \widetilde{H}_n(X^C;\bbZ)$ is bijective for $n \in \bbZ$.
Since $\widetilde{H}_n(X^C;\bbZ)$ is almost finitely generated,
it must be almost trivial for $n \ge 1$. This implies that
$\widetilde{H}_n(X^C;\bbQ) = 0$ for all $n \in \bbZ$.
\\[2mm]
\eqref{lem: Smith theory: almost isomorphism} This follows from
the Lerray-Serre spectral sequence of the fibration $X^C \to EC_GC
\times X^C \to BC_GC$ whose $E^2$-term is
$H_s(BC_GC;\widetilde{H}_t(X^C;\bbZ))$ and which converges to
$H_{s+t}(\pr \colon EC_GC \times_{C_GC} X^C \to BC_GC;\bbZ)$ and
the fact that the full subcategory of almost trivial
$\bbZ$-modules is a Serre subcategory of the abelian category of
$\bbZ$-modules.
%\\
%\eqref{lem: Smith theory: rational isomorphism} This is proved analogously to
%\eqref{lem: Smith theory: almost isomorphism}.
\end{proof}

Now we can give the proof of Theorem~\ref{the: Weakening the finiteness assumption}.
\begin{proof}
If one goes through the proofs of
Theorem~\ref{the: main theorem} and
Theorem~\ref{the: Multiplicative structure}
one sees that the finiteness condition about $\underline{E}G$ enters only, when we apply
Lemma~\ref{lem: i^*_G(X;bfE) is bijective for finite X} to $\underline{E}G$. Hence it suffices to show
that under the assumptions appearing in Theorem~\ref{the: Weakening the finiteness assumption}
the map
$$i^n_G(\underline{E}G;\bfK) \colon H_G^n\left(\underline{E}G;\bfK_{\Bor}\right) \otimes_{\bbZ} \bbQ \to
H_G^n\left(\underline{E}G;\left(\bfK_{\Bor}\right)_{(0)}\right)$$
is a $\bbQ$-isomorphism for all $n \in \bbZ$.

Let $f \colon X \to \underline{E}G$ be the up to $G$-homotopy unique $G$-map.
We obtain the following commutative diagram
\comsquare{H_G^n\left(\underline{E}G;\bfK_{\Bor}\right) \otimes_{\bbZ} \bbQ}
{i^n_G(\underline{E}G;\bfK)}
{H_G^n\left(\underline{E}G;\left(\bfK_{\Bor}\right)_{(0)}\right)}
{H_G^n\left(f;\bfK_{\Bor}\right) \otimes_{\bbZ} \bbQ}
{H_G^n\left(\underline{E}G;\left(\bfK_{\Bor}\right)_{(0)}\right)}
{H_G^n\left(X;\bfK_{\Bor}\right) \otimes_{\bbZ} \bbQ}
{i^n_G(X;\bfK)}
{H_G^n\left(X;\left(\bfK_{\Bor}\right)_{(0)}\right)}
Since $H_k(f,\bbZ) \colon H_k(X;\bbZ) \to H_k(\underline{E}G;\bbZ)$ is bijective for
all $k \in \bbZ$, we conclude from the Lerray-Serre spectral sequence that
$H_k(EG \times_G f,\bbZ) \colon H_k((EG \times_G X;\bbZ) \to H_k((EG \times_G\underline{E}G;\bbZ)$ is bijective for
all $k \in \bbZ$. This implies that the left vertical arrow in the commutative square above
which can be identified with $K^n(EG \times_G f) \colon K_n(EG \times_G \underline{E}G) \to
K_n(EG \times_G X)$ is bijective.
The lower horizontal arrow is bijective by Lemma~\ref{lem: i^*_G(X;bfK) is bijective for certain X}.
The right vertical arrow is bijective by Lemma~\ref{H^*,-,(K_Bor)_(0))-iso}
and Lemma~\ref{lem: Smith theory}~\eqref{lem: Smith theory: almost trivial}.
Hence the upper horizontal arrow is bijective.

The claim about the equivalent reformulation of condition
\eqref{the: Weakening the finiteness assumption: condition 4}
follows from
Lemma~\ref{lem: almost iso induced by EG times_G X to X/G}
\eqref{lem: almost iso induced by EG times_G X to X/G: almost isomorphism}
and
Lemma~\ref{lem: Smith theory}~\eqref{lem: Smith theory: almost isomorphism}.
This finishes the proof of
Theorem~\ref{the: Weakening the finiteness assumption}.
\end{proof}

%%%%%%%%%%%%%%%%%%%%%%%%%%% Section 7 %%%%%%%%%%%%%%%%%%%%%%%%%%%%%%%%%%
\typeout{-----------------------  Section 7 ------------------------}

\tit{Examples and Further Remarks}

Some finiteness conditions such as appearing in Theorem~\ref{the: Weakening the finiteness assumption}
are necessary as the following example shows.

\begin{example}{\bf (Necessity of the finiteness conditions).} \label{exa: some finiteness conditions are crucial.}
\em
 Consider
$G = \ast_{i=1}^{\infty} \bbZ/p$ for a prime number $p$.
Then $BG \simeq \bigvee_{i=1}^{\infty} B\bbZ/p$ and we get
$$K^0(BG) ~ \cong ~ K^0(\pt) \times \prod_{i=1}^{\infty} \widetilde{K}^0(B\bbZ/p) ~ \cong ~
\bbZ \times \prod_{i=1}^{\infty} (\bbZ\widehat{_p})^{p-1},$$
if $\pt$ is the one-point-space.
Since $H^n(BG;M) \cong \prod_{i=1}^{\infty} H^n(B\bbZ/p;M) = 0$ for any $\bbQ$-module $M$ and
$n \ge 2$, the cohomological dimension of $G$ over $\bbQ$ is $\le 1$ and hence
$G$ acts on a tree $T$ with finite stabilizers~\cite{Dunwoody(1979)}.
Then $T$ is a $1$-dimensional model for $\underline{E}G$
(see~\cite[page~20]{Serre(1980)} or
\cite[Proposition~4.7 on page~17]{Dicks-Dunwoody(1989)}). By the
Kurosh Subgroup Theorem
\cite[Theorem~1.10 on page~178]{Lyndon-Schupp(1977)})
any non-trivial finite subgroup of $G$ is conjugated to precisely one of the summands $\bbZ/p$ and is equal
to its centralizer. Hence
$p$ is an upper bound on the orders of finite subgroups of $G$. Obviously
$\bbZ\widehat{_p} \otimes_{\bbZ} \bbQ$ is canonically isomorphic to $\bbQ\widehat{_p}$.
If the conclusion of Theorem~\ref{the: main theorem} would be true for $G$, it would predict that the canonical map
$$\left(\prod_{i=1}^{\infty} (\bbZ\widehat{_p})^{p-1}\right) \otimes_{\bbZ} \bbQ ~ \to ~
\prod_{i=1}^{\infty} \left((\bbZ\widehat{_p})^{p-1} \otimes_{\bbZ} \bbQ\right) ~ = ~
\prod_{i=1}^{\infty} (\bbQ\widehat{_p})^{p-1}$$
is bijective, what is not true. For instance, the element $(p^{-i})_{i=1}^{\infty}$ is not contained in its image.

Notice that in this example all conditions appearing in
Theorem~\ref{the: Weakening the finiteness assumption} are
satisfied except the condition that the set of conjugacy classes
$(C)$ of finite cyclic subgroups $C \subseteq G$ of prime power
order with $T^C \not= \emptyset$ is finite;

We emphasize that no restriction (except properness) occur in
Lemma~\ref{lem: refined non-multiplicative rational computation of H^*_G(X;(bfK_{Bor})_{(0)})}.
The problem is in Lemma~\ref{H^*,-,(K_Bor)_(0))-iso} some additional finiteness assumptions are needed.
\em
\end{example}

\begin{example}[$SL_3(\bbZ)$] \label{exa: SL_3(Z)} \em
Consider the group $G = SL_3(\bbZ)$.
It is well-known that its rational cohomology satisfies
$\widetilde{H}^n(BSL_3(\bbZ);\bbQ) = 0$ for all $n \in \bbZ$.
Actually, we conclude from ~\cite[Corollary on page~8]{Soule(1978)}
that for $G = SL_3(\bbZ)$ the quotient space $G\backslash\underline{E}G$ is contractible
and compact.
From the classification of finite subgroups of $SL_3(\bbZ)$ we
see that $SL_3(\bbZ)$ contains up to conjugacy  two elements of order $2$,
two elements of order $4$ and two elements of order $3$ and no further conjugacy classes of
non-trivial elements of prime power order.  The rational homology of all the
centralizers of elements in $\con_2(G)$ and $\con_3(G)$ agree with the one of the trivial group
(see~\cite[Example~6.6]{Adem(1993b)}). Hence Theorem~\ref{the: main theorem} shows
\begin{eqnarray*}
K^0(BSL_3(\bbZ)) \otimes_{\bbZ} \bbQ & \cong &
\bbQ \times (\bbQ\widehat{_2})^4 \times (\bbQ\widehat{_3})^2;
\\
K^1(BSL_3(\bbZ)) \otimes_{\bbZ} \bbQ & \cong &  0.
\end{eqnarray*}
The identification of $K^0(BSL_3(\bbZ)) \otimes_{\bbZ} \bbQ$ above
is compatible with the multiplicative structure on the target
described in Example~\ref{exa: including multiplicative structures}.

Actually the computation using Brown-Petersen cohomology
and the Conner-Floyd relation by Tezuka and Yagita
\cite{Tezuka-Yagita(1992)} gives the integral computation
\begin{eqnarray*}
K^0(BSL_3(\bbZ)) & \cong &
\bbZ \times (\bbZ\widehat{_2})^4 \times (\bbZ\widehat{_3})^2;
\\
K^1(BSL_3(\bbZ)) & \cong &  0.
\end{eqnarray*}
\em
\end{example}

\begin{example}[Groups with appropriate maximal finite subgroups] \label{exa: conditions M and NM} \em
Let $G$ be a discrete group.
Consider the following assertions concerning $G$:
\begin{itemize}

\item[(M)] Every non-trivial finite subgroup of $G$ is contained in a unique maximal finite subgroup;

\item[ (NM)] If $M \subseteq G$ is maximal finite, then $N_GM = M$;

\item[(C)] There is a cocompact model for $\underline{E}G$.
\end{itemize}

The conditions (M) and (NM) imply the following: For any  non-trivial finite subgroup $H \subseteq G$
we have $N_GH = N_MH$ if $M$ is a maximal finite subgroup containing $H$.
Let $\{M_i \mid i \in I\}$ be a complete set of representatives of the conjugacy classes of maximal finite
subgroups of $G$. Fix a prime $p$. Then the obvious map
$$\coprod_{i \in I} \con_p(M_i) \xrightarrow{\cong} \con_p(G)$$
is a bijection. Let $r_p(M_i) = |\con_p(M_i)|$ be the number of
conjugacy classes of elements in $M_i$ of order $p^k$ for some $k \ge 1$.

Theorem~\ref{the: main theorem} yields for a group
satisfying conditions (M), (NM) and (C) above
rational isomorphisms
\begin{eqnarray*}
K^0(BG) \otimes_{\bbZ} \bbQ & \cong &
\prod_{i \in \bbZ} H^{2i}(BG;\bbQ) \times
\prod_{p} \left(\bbQ\widehat{_p}\right)^{\sum_{i \in I} r_p(M_i)}
\\
K^1(BG) \otimes_{\bbZ} \bbQ & \cong &
\prod_{i \in \bbZ} H^{2i+1}(BG;\bbQ).
\end{eqnarray*}

Here are some examples of groups $Q$ which satisfy conditions (M), (NM) and (C)
\begin{itemize}

\item Extensions $1 \to \bbZ^n \to G \to F \to 1$ for finite $F$ such that the conjugation
  action of $F$ on $\bbZ^n$ is free outside $0 \in \bbZ^n$. \\[1mm]
The conditions (M), (NM)  are satisfied by
\cite[Lemma~6.3]{Lueck-Stamm(2000)}. There are models for $\underline{E}G$
whose underlying space is $\bbR^n$. The quotient $G\backslash\underline{E}G$ looks like the quotient
of $T^n$ by a finite group.

\item Fuchsian groups $F$ \\[1mm]
See for instance~\cite[Lemma~4.5]{Lueck-Stamm(2000)}).
The quotients $G\backslash \underline{E}G$ are closed orientable surfaces.
In ~\cite{Lueck-Stamm(2000)} the larger class of cocompact planar
groups (sometimes also called cocompact NEC-groups) is treated.

\item Finitely generated one-relator groups $G$\\[1mm]
Let $G = \langle (q_i)_{i \in I} \mid r \rangle$
be a presentation with one relation.  We only
have to consider the case, where $Q$ contains torsion.
Let $F$ be the free group with basis $\{q_i \mid i \in I\}$. Then $r$ is an element in
$F$. There exists an element $s \in F$ and an integer $m \ge 2$
such that $r = s^m$, the cyclic  subgroup $C$
generated by the class $\overline{s} \in Q$ represented
by $s$ has order $m$, any finite subgroup of $G$ is subconjugated to $C$ and
for any $q \in Q$ the implication
$q^{-1}Cq \cap C \not= 1 \Rightarrow q \in C$ holds.
These claims follows from
\cite[Propositions~5.17, 5.18 and~5.19 in~II.5 on pages~107 and~108]{Lyndon-Schupp(1977)}.
Hence $Q$ satisfies (M) and (NM). There are explicit two-dimensional models
for $\underline{E}G$ with one $0$-cell $G/C \times D^0 $, as many free $1$-cells $G \times D^1$  as there are elements in $I$ and
one free $2$-cell $G \times D^2$ (see~\cite[Exercise~2~(c)~II. 5 on page~44]{Brown(1982)}).
\end{itemize}

For the three examples above one can make
$H^*(BG;\bbQ) = H^*(G\backslash\underline{E}G;\bbQ)$ more explicit.
\em
\end{example}

\begin{example}[Extensions of $\bbZ^n$ with $\bbZ/p$ as quotient] \label{exa: extensions by Z/p} \em
Suppose that $G$ can be written as an extension
$1 \to A \to G \to \bbZ/p \to 1$ for some fixed prime number $p$
and for $A = \bbZ^n$ for some integer $n \ge 0$ and that $G$ is not torsionfree.
The conjugation action of $G$ on the normal subgroup $A$ yields
the structure of a $\bbZ[\bbZ/p]$-module on $A$. Every
non-trivial element $g \in G$ of finite order $G$ has order $p$ and satisfies
$$N_G\langle g  \rangle = C_G\langle g \rangle = A^{\bbZ/p} \times \langle g \rangle.$$
There is a bijection
$$\mu \colon H^1(\bbZ/p;A) \times (\bbZ/p)^\times ~ \xrightarrow{\cong} ~ \con_p(G),$$
where $H^1(\bbZ/p;A)$ is the first cohomology of $\bbZ/p$ with coefficients
in the $\bbZ[\bbZ/p]$-module $A$. If we fix an element $g \in G$ of order $p$ and a generator $s \in \bbZ/p$,
the bijection $\mu$  sends $([u],\overline{k}) \in  H^1(\bbZ/p;A) \times (\bbZ/p)^\times$
to the conjugacy class $(ug^k)$ of $ag^k$ if $[u] \in H^1(\bbZ/p;A)$ is represented by the element
$u$ in the kernel of the second differential $A \to A, ~ a \mapsto \sum_{i=0}^{p-1} s^i \cdot a$
and $k \in \bbZ$ represents $\overline{k}$. There is a cocompact model for $\underline{E}G$ with $A \otimes_{\bbZ} \bbR$
as underlying space.
Hence Theorem~\ref{the: main theorem} yields for $G$ as above
rational isomorphisms
$$
K^n(BG) \otimes_{\bbZ} \bbQ ~ \cong ~
\prod_{i \in \bbZ} H^{2i+n}(BA;\bbQ)^{\bbZ/p} \times
\prod_{k=1}^{r} ~ \prod_{i \in \bbZ} H^{2i+n}\left(B(A^{\bbZ/p});\bbQ\widehat{_p}\right),
$$
if we put $r = (p-1) \cdot |H^1(\bbZ/p;A)|$.

Take for instance $A$ to be the cokernel of the inclusion of
$\bbZ[\bbZ/p]$-modules $\bbZ \to \bbZ[\bbZ/p], n \mapsto n \cdot \sum_{i=0}^{p-1} t^i$,
where $\bbZ$ carries the trivial $\bbZ/p$-action and $t \in \bbZ/p$ is a fixed generator.
One can identify $A$ with the extension of $\bbZ$ by adjoining a primitive $p$-th root of unity.
From the long exact cohomology sequence associated to the short exact sequence of
$\bbZ[\bbZ/p]$-modules $0 \to \bbZ \to \bbZ[\bbZ/p] \to A \to 0$ one concludes
that $H^1(\bbZ/p;A)$ is a cyclic group of order $p$. One easily checks that
$A^{\bbZ/p} = 0$. Hence we obtain for the semi-direct product $A \rtimes \bbZ/p$
\begin{eqnarray*}
K^0(B(A \rtimes \bbZ/p)) \otimes_{\bbZ} \bbQ &  \cong &
\bbQ \times (\bbQ\widehat{_p})^{p^2-p};
\\
K^1(B(A \rtimes \bbZ/p)) \otimes_{\bbZ} \bbQ &  \cong & 0.
\end{eqnarray*}
The identification of $K^0(B(A \rtimes \bbZ/p)) \otimes_{\bbZ} \bbQ$
is compatible with the multiplicative structure on the target
described in Example~\ref{exa: including multiplicative structures}.
\em
\end{example}

\begin{remark}[Comparison with Adem's work]  \label{rem: comparison with Adem} \em
The results and the examples appearing in this paper are consistent
with the ones by Adem ~\cite{Adem(1993b)}. Adem needs that $G$ contains a normal torsionfree subgroup
$G' \subseteq G$ of finite index and uses the Atiyah-Segal Completion Theorem for the finite group $G/G'$
to compute rationally the $K$-theory with coefficients in the $p$-adic integers $\bbZ\widehat{_p}$.
His condition that in his notation $\Gamma'\backslash X$ is compact is precisely the condition
that there is a cocompact model for $\underline{E}G$. We can drop the condition
of the existence of the normal torsionfree subgroup $G' \subseteq G$ of finite index with our methods.

One can get Adem's  local computations from ours by replacing for a fixed prime
$p$ the cohomology $H^*(BG;\bbQ)$ by $H^*(BG;\bbQ\widehat{_p})$ and ignoring in the product
running over all primes all the contributions coming from primes different from $p$. For instance
Example ~\ref{exa: SL_3(Z)} implies that the $\bbQ\widehat{_3}$-algebra
$K^0(BSL(3,\bbZ);\bbZ\widehat{_3}) \otimes_{\bbZ\widehat{_3}} \bbQ\widehat{_3}$
is given by
$$\bbQ\widehat{_3}[u,v]/(u^2 = u, v^2= v, uv = 0).$$
If one makes the change of variables $u = (\alpha - 2)/3$ and $v = (\beta - 2)/3$, one obtains the
presentation in~\cite[Example~6.6]{Adem(1993b)}
$$\bbQ\widehat{_3}[u,v]/(\alpha^2 - \alpha -2, \beta^2 -\beta - 2, \alpha\beta-2(\alpha + \beta -2)).$$
Recall that after complexification we can determine the multiplicative structure in general
(see Theorem~\ref{the: Multiplicative structure}).

There are interesting discussions about Euler characteristics and maps of groups inducing isomorphisms
on homology in~\cite{Adem(1992)} and ~\cite{Adem(1993b)} which also apply to our setting.
\em
\end{remark}

\begin{remark}[Hodgkin's computation]
\label{rem: Hodgkins computation} \em
Let $\Gamma^n$ be the mapping class group of the sphere $S^2$ with $n$ punctures for $n \ge 3$.
Hodgkin computes rationally the $K$-theory of $B\Gamma^n$ with coefficients in the $p$-adic integers $\bbZ\widehat{_p}$
using Adem's formula in~\cite{Adem(1992)}. The main work done in the paper by
Hodgkin~\cite[Proposition~2.2 and Theorem~2]{Hodgkin(1995)} is to figure out
the set of conjugacy classes of elements of order $p^s$ for each prime $p$ and integer $s \ge 1$
and the rank of $K^k(BC_G \langle g \rangle) \otimes_{\bbZ} \bbQ \cong
\prod_{i \in \bbZ} H^{2i + k}(BC_G\langle g \rangle;\bbQ)$ for each element $g \in \Gamma^n$
of prime power order. One can identify
$K^*(BC_G \langle g \rangle) \otimes_{\bbZ} \bbQ$ with $(K^*(BK^r) \otimes_{\bbZ} \bbQ)^{\Sigma_r}$
or  with $(K^*(BK^r) \otimes_{\bbZ} \bbQ)^{\Sigma_{r-2} \times \Sigma_2}$ for appropriate
integers $r$ depending only on the order of $g$, where
$K^r$ is the pure mapping class group of $S^2$ with $r$ punctures and $\Sigma_t$ denotes
the symmetric group of permutation of the set consisting of $t$ elements.
Thus one obtains with Theorem~\ref{the: main theorem} the precise structure
of the $\bbQ$-vector spaces $K^k(B\Gamma^n) \otimes_{\bbZ} \bbQ$. It may be worthwhile
to investigate the product structure on the cohomology $H^*(C_G\langle g \rangle;\bbC)$ since this
would lead to a computation of $K^*(B\Gamma^n) \otimes_{\bbZ} \bbC$ including its multiplicative structure by
Theorem~\ref{the: Multiplicative structure}.
\em
\end{remark}

\begin{remark}[Criterion for torsionfree] \label{rem: characterization of torsionfree}
\em
Let $G$ be a discrete group with a finite model for $\underline{E}G$.
Then the following assertions are equivalent:
\begin{enumerate}
\item $G$ is torsionfree;

\item The abelian group $K^k(BG)$ is finitely generated for $k \in \bbZ$;

\item The rational vector space $K^k(BG) \otimes_{\bbZ} \bbQ$ is finite dimensional for $k \in \bbZ$.

\end{enumerate}
An application of the Atiyah-Hirzebruch spectral sequence proves the implication
(a) $\Rightarrow$ (b). The implication (b) $\Rightarrow $(c) is obvious. The implication (c) $\Rightarrow $(a) follows from
Theorem~\ref{the: main theorem}  since $\bbQ\widehat{_p}$ is an infinite dimensional $\bbQ$-vector space.
\em
\end{remark}

\begin{remark}[Torsion prime to the order of finite subgroups] \em
Suppose that there is a finite model for $\underline{E}G$. Let the
ring $\Lambda^G$ be the subring of $\bbQ$ obtained by inverting
the orders of finite subgroups of $G$. We state without giving the
details of the proof but referring to~\cite{Joachim-Lueck(2005)}
that one can improve Theorem~\ref{the: main theorem} to the
statement that there is a $\Lambda^G$-isomorphism
\begin{multline*}
\overline{\ch}^n_{G,\lambda^G} \colon K^n(BG) \otimes_{\bbZ} \Lambda^G
~  \xrightarrow{\cong}
\\
K^n(G\backslash \underline{E}G) \otimes_{\bbZ} \Lambda^G \times
\prod_{p \text{ prime}} ~ \prod_{(g) \in \con_p(G)}
\left(\prod_{i \in \bbZ} H^{2i+n}(BC_G\langle g \rangle;\bbQ\widehat{_p})\right).
\end{multline*}
Consider a prime $q$ for which there exists no element order $q^s$
for some $s \ge 1$ in $G$, in other words, $q$ is not invertible
in $\Lambda^G$. Then the projection $BG = EG \times_G
\underline{E}G \to G\backslash \underline{E}G$ induces an
isomorphism
$$\tors_q\left(K^n(G\backslash \underline{E}G)\right) ~ \xrightarrow{\cong} \tors_q\left(K^n(BG)\right),$$
if for an abelian group $A$ we denote by $\tors_q(A)$ the subgroup of elements
$a \in A$ which are annihilated by some power of $q$. In particular
$K^n(BG)$ contains $q$-torsion if and only if
$K^n(G\backslash \underline{E}G)$  contains $q$-torsion.
It can occur that $K^n(G\backslash \underline{E}G)$
contains elements of order $q^s$ for some $s \ge 1$
(see~\cite{Leary-Nucinkis(2001a)}.)
We will explain in~\cite{Joachim-Lueck(2005)} that
the subgroup of torsion elements in $K^n(BG)$ is finite.
\em
\end{remark}

%%%%%%%%%%%%%%%%%%%%%%%%%%%%%%%%%%%%%%%  References %%%%%%%%%%%%%%%%%%%%%%%%%%%%%%%%%%

\bibliographystyle{abbrv}
\bibliography{dbdef,dbpub,dbpre,dbtkcsrextra}

\end{document}